# End-to-End Integrated Simulation for Predicting the Fate of Contaminant and Remediating Nano-Particles in a Polluted Aquifer

Shikhar Nilabh, Fidel Grandia

Amphos 21 Consulting S.L., c/Veneçuela 103, 08019 Barcelona, Spain

## Abstract
Groundwater contamination caused by Dense Non-Aqueous Phase Liquid (DNAPL) has an adverse impact on human health and environment. Remediation techniques, such as the in-situ injection of nano Zero Valent Iron (nZVI) particles, are widely used in mitigating DNAPL-induced groundwater contamination. However, an effective remediation strategy requires predictive insights and understanding of the physiochemical interaction of nZVI and contamination along with the porous media properties. While several stand-alone models are widely used for predictive modeling, the integration of these models for better scalability and accuracy is still rarely utilized. This study presents an end-to-end integrated modeling framework for the remediation of DNAPL-contaminated aquifers using nZVI. The framework simulates the migration pathway of DNAPL and subsequently its dissolution in groundwater resulting in an aqueous contaminant plume. Additionally, the framework includes simulations of nZVI mobility, transport, and reactive behavior, allowing for the prediction of the radius of influence and efficiency of nZVI for contaminant degradation. The framework has been applied to a hypothetical 2-dimensional and heterogeneous silty sand aquifer, considering trichloroethylene (TCE) as the DNAPL contaminant and carboxymethyl cellulose (CMC) coated nZVI for remediation. The results demonstrate the framework's capability to provide comprehensive insights into the contaminant's behavior and the effectiveness of the remediation strategy. The proposed modeling framework serves as a reference for future studies and can be expanded to incorporate real field data and complex geometries for upscaled applications.

## 1 Introduction
Groundwater contamination from industrial releases or discharges poses a serious threat to human health and the environment. An important form of this pollution is dense non-aqueous phase liquid (DNAPL) that infiltrates into the subsurface with a tendency to migrate downward. The DNAPL phase immobilizes as pool and ganglia resulting in long-term contamination of groundwater due to their dissolution. The elimination of these contaminants from the groundwater often requires active remediation techniques. Particularly, in situ injection of nano zero valent iron (nZVI) particles has proven to be a successful strategy for remediation of groundwater contamination in previous studies (Sethi and Di Molfetta, 2019). These particles are suitable for groundwater remediation due to their rapid reaction kinetics, minimal ecological impact, nontoxicity, and cost-effectiveness (Karn et al.,2009, Zhao et al.,2016). Owing to its high efficiency, nZVI has been widely used for the remediation of aquifers contaminated with DNAPLs such as organohalides and mercury (D'Aniello, 2017).

Designing a remediation plan using nZVI requires an efficient contaminant characterization that includes identifying the extent of contamination zones and understanding the governing factors influencing their dynamic behavior. Modeling is often used to identify the dynamic contamination zones within the aquifer and optimize the design of remediation systems, thereby minimizing the risk of failure. For contaminant characterization, simulations based on two-phase flow, contamination dissolution transport model are widely used. Two-phase models are instrumental in computational reconstruction of the migration pathways and prediction of DNAPL source zone architecture (Helmig *et al.*, 2006; Kamon *et al.*, 2004; Kueper and Frind, 1991; McLaren *et al.*, 2012; Zheng *et al.*, 2015). Furthermore, several models have been built to predict the dissolution of DNAPL source zone and aqueous contamination plume geometry (Almpanis *et al.*, 2021; Hwang *et al.*, 2013; Kokkinaki *et al.*, 2013). The prediction from these models gives useful insight for the contamination longevity, aqueous DNAPL plume geometry, and their spatio-temporal evolution.

On one hand, several researchers have developed computational tools to study DNAPL migration, its dissolution, and transport. Extensive models have also been built to characterize the mobility and reactivity of the nZVI in field and laboratory scale. (Babakhani et al., 2018; Krol et al., 2013; Tosco et al., 2014; Xu et al., 2019). Numerical simulation of ZVI fate is used to obtain predictive insights into the role of physio-chemical mechanisms governing the nZVI mobility which is useful for remediation design. In addition, these simulations estimate the Radius of Influence of nZVI for contaminant degradation as well as its efficiency and timeline for groundwater remediation. While these numerical simulators are useful as stand-alone models in characterizing contamination and nZVI behavior enabling its extensive application, less emphasis has been given to the integrated modeling framework. An end-to-end integrated modeling framework can provide significant benefits in designing a remediation plan. First, it can lead to greater accuracy by accounting for complex interactions and feedback loops between different phases of remediation. Second, it provides better scalability because it can handle different scales of laboratory and field data. Third, it improves interpretability by providing a more complete and coherent picture of the field contamination and remediation plan. Finally, it provides a common modeling platform for different researchers and experts working on different components of remediation. This allows them to share data, findings, and feedback in a more structured and coherent way.

Few attempts have been made in previous studies to model the different phases of the groundwater remediation study. For example, Kokkinaki (2013) developed an integrated model for the infiltration and rate limited dissolution of DNAPL. The model was extended to simulate the dissolved component of DNAPL geometry governed by advection and dispersion processes in the groundwater. Likewise, Hwang et al. (2013) extended the DNAPL migration model to simulate the reactive and transport of the aqueous DNAPL origination from the DNAPL dissolution. Fagerlund et al. (2012) developed an integrated column scale model for simulating the dissolution transport and degradation of DNAPL in the presence of nZVI. Tsakiroglou et al. (2018) integrated the simulation of ZVI transport As well as PCE dissolution, transport and degradation in a column scale sand. Taghavy et al. (2010) developed a column scale multiphase transport model that accounted for rate-limited dissolution of DNAPL, and its reaction with the nZVI deposited on the sand grain surface. While these models have been successful in integrating different stages of remediation and planning a suitable remediation design, an end-to-end integration of these models is still missing.

To our knowledge, this paper represents the first documented attempt to develop an end-to-end continuum-scale modeling framework that provides relevant insights for all relevant stages of remediation strategy, from the entry of DNAPL into the aquifer to the degradation of the contaminant. In the present work, a hypothetical 2-dimensional sandy aquifer of 35 x 12 m dimension is considered. The domain includes micro and macro-scale heterogeneities to simulate the real field-like condition. The model considers DNAPL in the form of

Trichloroethylene (TCE), a volatile organic compound which is a common form of industrial pollutant found in contaminated groundwater. Several researchers have reported an enhanced transport of ZVI with surface coating of CMC and has been widely used for groundwater remediation. Hence, CMC coated nZVI has been considered for the contaminant degradation in the hypothetical aquifer. The model has been intended to cover all the major stages encountered when assessing a contaminated site with the objective of planning a pilot remediation operation. Few reasonable assumptions for the parameter value have also been made to characterize the fate of TCE and ZVI in the hypothetical 2D aquifer. Following are the stages that have been simulated:

(1) Stage 1 deals with the simulation of the TCE phase migration pathway and the formation of source zone architecture. Given the sparingly soluble nature of TCE, its dissolution during infiltration in the aquifer is assumed to be negligible. The model further assumes negligible TCE mass transfer due to adsorption or volatility. The simulation results provide the morphology and quantification of the TCE ganglia and pool, as well as a prediction of the source zone architecture.
(2) Stage 2, the dissolution of TCE and the transport of the aqueous TCE is simulated to predict the dynamic contaminant plume geometry in the groundwater flow. Under the integrated modeling framework, the TCE source zone architecture predicted in Stage 1 is used as the initial condition for the dissolution and transport simulation. This simulation stage also accounts for the evolution of TCE source zone with the dissolution, highlighting its gradual depletion over time.
(3) Stage 3 of the study involves the simulation of a pilot-scale injection using a CMC-ZVI solution at the downgradient of a TCE pool predicted in Stage 1. This simulation is designed to facilitate a deeper understanding of the transport and retention mechanisms of nano zero-valent iron (nZVI) within the hypothetical aquifer. The model result quantifies the Radius of Influence (ROI) for ZVI, including its spatial distribution within the aquifer. Furthermore, the clogging effect and viscosity change due to CMC-nZVI injection are simulated in a fully coupled approach.
(4) Stage 4 of the study focuses on reactive transport modeling, which simulates the interaction between the aqueous contaminant and nZVI leading to the contaminant degradation. The resulting model provides insight into the dynamics of the aqueous TCE plume geometry in the presence of nZVI. Furthermore, this stage of the study also simulates the depletion of active nZVI surface for TCE degradation, as well as the simulation of the TCE source zone.

Therefore, the exercise aims to demonstrate an end-to-end framework that encapsulates all the major stages of remediation. The model can serve as a reference platform where the parameters or simulation can be extended according to the needs of the model. The numerical model in this work is built on the Comsol Multiphysics software package (comsol 2021), which uses Galerkin based Finite Element method for solving Partial Differential Equations (PDEs). The software´s capability of coupling diverse physical phenomena as well as customizing the mesh and solver for each problem make it suitable for developing the robust two-phase flow model. The software uses MUltifrontal Massively Parallel sparse (MUMPs) Solver, a direct method for solving the linear system of equations obtained from Finite Element Analysis. The model represents the flow dynamics at a macroscopic scale which allows continuum-based parameters and variables in Comsol.

## 2 Implementation of an integrated modeling framework

### 2.1 Geometry and sedimentology of the polluted aquifer

The integrated model domain consists of a 2D, unconfined aquifer with dimensions of 35 meters width and 12 meters height (Figure 1). To incorporate the natural complexity found in a common aquifer, macroscopic heterogeneity has been introduced by two layers of sand with different hydrogeological properties and 6 clay lenses. The aquifer computationally represents two silty sand layers and six clay lenses assigned at random location within the model domain. The average permeability of $1\times10^{-12}$ m$^2$ for the upper silty sand layer and $0.5\times10^{-12}$ m$^2$ for lower silty sand layer. A constant permeability of $5\times10^{-14}$ m$^2$ is assumed for clay. These values are taken considering the usual permeability and porosity range of the silty sand and clays (Krol *et al.*, 2013; Lewis *et al.*, 2006). An arbitrary but reasonable value of entry pressure, residual saturation of water and TCE have also been assigned for the silty sand and clay. These values lie within the range of the previously observed values for silty sand and clay (Kueper and Frind, 1991). Figure 1 shows the hydrogeological properties used for the simulation of contaminant´s fate within the aquifer.

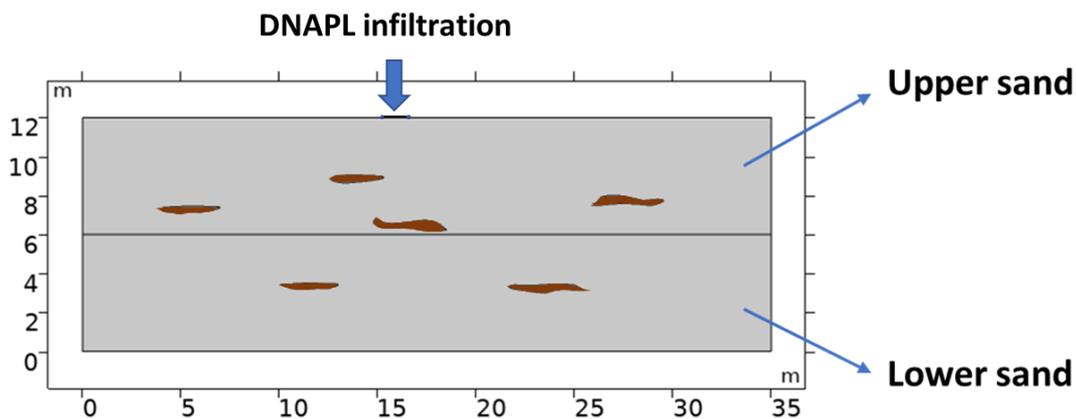

*Figure 1: Model domain for the simulation of the integrated stages of contamination and remediation. The main components are two sand layers and six clay units (illustrated with brown color) arranged randomly to represent a typical field case.*

Microscale heterogeneity also influences the flow and transport of contamination. Therefore, to study the impact of microscale heterogeneity on contamination and remediation dynamics, a statistical permeability field is generated (Figure 2). The field was generated using the Field Generator in PMWIN (Processing Modflow 5.3) that uses Mejia's algorithm (Frenzel, 1995; Mejía and Rodríguez-Iturbe, 1974). Micro-scale heterogeneity is assumed to have a log-normal distribution with log-variance of 0.2 for both the sand layers. The generated permeability field was imported into COMSOL utilizing the interpolation function, which facilitated direct data input without altering the discretization. No microscopic heterogeneity is considered in the clay lenses and thus a constant permeability is assigned to it.

*Table 1: Averaged hydrogeological properties of different geological units (upper sand layer, lower sand layer and clay lenses) in the aquifer used for the integrated modelling.*

| Properties | Upper layer | Lower Layer | Clay lenses |
|---|---|---|---|
| Mean permeability (m$^2$) | 1×10$^{-12}$ | 0.5×10$^{-12}$ | 5×10$^{-14}$ |
| Porosity | 0.4 | 0.27 | 0.25 |
| Mean residual saturation of water | 0.08 | 0.04 | 0.189 |
| Mean residual saturation of TCE | 0.08 | 0.04 | 0.04 |
| Entry pressure (Pa) | 1300 | 1500 | 3200 |

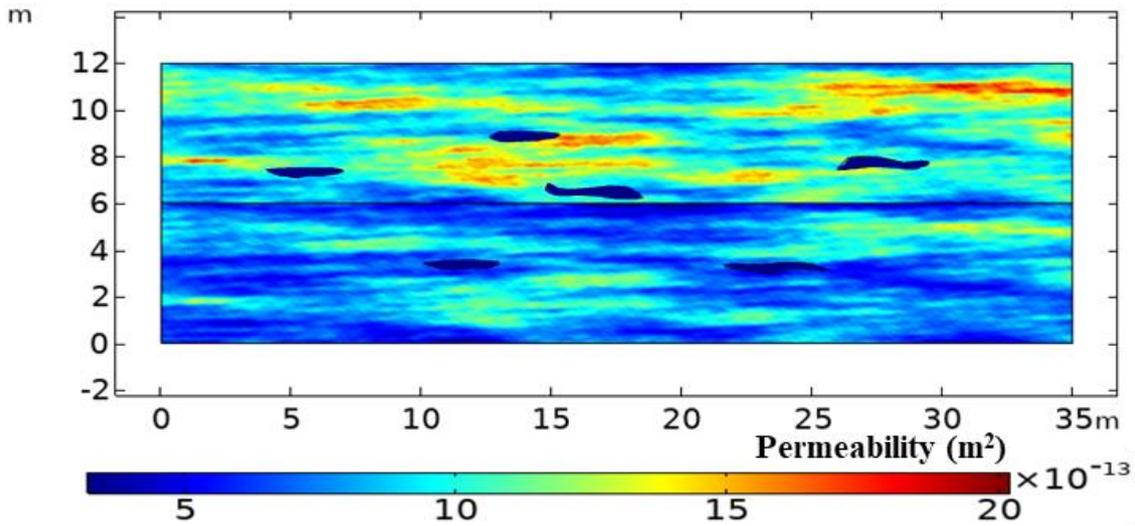

*Figure 2: Permeability distribution of the porous media in the integrated model, calculated from a logarithmic distribution.*

## 2.2 Contaminant properties

Trichloroethylene (TCE) is known to be denser than water and less viscous, thus it is expected to migrate vertically downward in the subsurface. Both water and TCE were assumed to be incompressible, with their densities held constant throughout the simulation. Previous studies have reported the TCE solubility in water to be 1.27 kg$_{TCE}$·m$^{-3}_{water}$ under natural groundwater condition (Grathwohl, 2012). Due to the low solubility of TCE in water, the dissolution of TCE in groundwater during the infiltration of the free phase was considered negligible during the initial modeling stage (Stage 1). The chemical and physical properties of both water and TCE are summarized in Table 1.

*Table 2: Fluid properties used in the simulation of TCE infiltration and dissolution in the groundwater (Grathwohl, 2012).*

| Fluid properties | value |
|---|---|
| Water density | 1000 kg·m$^{-3}$ |
| TCE density | 1470 kg·m$^{-3}$ |
| Water viscosity | 0.001 Pa·s |
| TCE viscosity | 0.0005 Pa·s |
| TCE dissolution in water | 1.27 kg·m$^{-3}$ |

## 2.3 Modelling Stage 1: Two-phase flow during contaminant infiltration

In Stage 1, the entry of TCE phase, its migration and immobilization in the hypothetical aquifer is simulated. The simultaneous two-phase flow equation governs the flow of water and TCE in the hypothetical aquifer (Equation 1,2). The Brooks and Corey relationship for capillary pressure and saturation has been used to close the system of equations. (Equation 3,4).

$$\theta \frac{\partial (S_w)}{\partial t} = \nabla \cdot \left( k \frac{k_{rw}}{\mu_w} (\nabla p_w - \rho_w g) \right) + q_w \qquad \text{(eq. 1)}$$

$$\theta \frac{\partial (S_n)}{\partial t} = \nabla \cdot \left( k \frac{k_{rn}}{\mu_n} (\nabla p_n - \rho_n g) \right) + q_n \qquad \text{(eq. 2)}$$

$$S_n + S_w = 1 \qquad \text{(eq. 3)}$$

$$p_c(S_w) = p_d S_e^{\frac{-1}{\lambda}} \qquad \text{(eq. 4)}$$

Where $S_n$ and $S_w$ are the saturation of TCE and water respectively, $p_n$ and $p_w$ are the pressure of TCE and water, $\theta$ is porosity, $q_w$ and $q_n$ are the volumetric source term for water and TCE respectively, $k$ is the permeability of the sand, $p_c$ is the capillary pressure defined as a function of water saturation. The TCE mass loss due to dissolution in groundwater, adsorption in the mineral rock, or natural degradation is assumed to be negligible, and therefore $q_w$ and $q_n$ are considered to be zero.

The numerical formulation for the above equations is developed on COMSOL. The literature reports different numerical schemes used for the simulation of two-phase flow in a porous media (Chen et al., 2006; Nilabh, 2021). The Implicit Phase Pressure and Explicit Saturation has been selected for numerical scheme due to its flexibility and decoupling in the PDE equation (Kueper and Frind, 1991). To incorporate the effect of clay lenses in the flow of TCE, Phase Pressure Saturation with Interface Condition (PPSIC) algorithm has been used The algorithm allows defining an explicit entry pressure at different lithological layers and has been used extensively in previous studies to predict the DNAPL accumulation at the lithological interfaces (Helmig et al., 2006).

### 2.3.1 Initial and boundary conditions

During the infiltration of TCE, the effect of groundwater flow on TCE migration path is considered to be negligible owing to the slow groundwater flow system in the hypothetical aquifer (see section 2.4.1). This is in accordance with the previous studies in the literature which reported minor effect of pore water on the DNAPL flow path for groundwater flow of range up to 10 m·s$^{-1}$ (Zheng *et al.*, 2015). Therefore, a hydrostatic initial condition with water table at the surface (hydraulic head =12 m) has been assumed. A Dirichlet boundary condition of 12 meters of hydraulic head is assigned at the left and right boundaries. The top and bottom boundaries are assigned with no flow boundary for both water and TCE, except in the area of infiltration.

The model domain is considered to be fully saturated with water at the initial stage of the simulation. The flow rate of the TCE phase from the top and the timeline has been chosen such that a significant amount of TCE phase and ganglia could be formed in the hypothetical aquifer. For this, the model has considered that TCE phase has infiltrated from the specified area at top of the aquifer (Figure 1) with an arbitrary flow rate of 0.001 kg·m$^{-2}$s$^{-1}$ for 35 days. After 35 days, the model simulates the further TCE phase migration within the aquifer which is governed by the fluid density difference and capillary force at the water-TCE phase interface. In total, Stage 1 simulates the TCE migration for 135 days after which the simulated TCE phase stabilizes as residual saturation along the migratory pathway and TCE pool at the top of impervious layers. The significant output from this stage is to predict the TCE source zone architecture which would be useful in subsequent stages for contaminant characterization and remediation designing.

## 2.4 Modelling Stage 2: TCE dissolution and transport

The TCE source zone in a real field case is expected to undergo slow depletion due to mass transfer of the contaminant to the aqueous phase. A dynamic plume of dissolved TCE is formed in the contaminated groundwater due to advection and dispersion processes. Simulating the dissolution and transport of TCE facilitates the understanding of the spatio-temporal dynamics of contamination in the groundwater hence aiding in an optimal remediation design. Therefore, in Stage 2 the model simulates the dissolution of TCE phase in water as well as the transport as dissolved phase in the hypothetical aquifer. The source zone architecture after 100 days of evolution predicted in Stage 1 has been taken as a starting point for Stage 2. The depletion of source zone due to dissolution is also simulated with TCE transport within the integrated framework. The model uses Advection Dispersion equation for simulating the transport of aqueous TCE in the groundwater (Equation 5):

$$\frac{\partial(\theta c)}{\partial t} = -\nabla(vc) + \nabla((D_d + \alpha v)\nabla c) + q \qquad \text{(eq. 5)}$$

Where $D_d$ is the molecular diffusion coefficient (m$^2$·s$^{-1}$), $v$ is the Darcy velocity (m·s$^{-1}$), and $c$ is the contaminant concentration (Kg·m$^{-3}$). $\alpha$ is the dispersivity(m), q is the mass flux (Kg·m$^{-2}$·s$^{-1}$), $K_l$ is the lumped mass transfer coefficient (m·s$^{-1}$), $C_s$ is the equilibrium concentration of contaminant in water (kg$_{TCE}$·m$^{-3}$).

Simultaneously, the depletion of TCE source zone in the hypothetical aquifer due to the dissolution is simulated using Equation 6:

$$\theta \frac{\partial(\rho_o S_n)}{\partial t} = q \qquad \text{(eq. 6)}$$

In some of the previous studies, Local Equilibrium assumption for contaminant dissolution has been observed valid in some field study where the aqueous DNAPL dissolves with solubility limit(Frind et al., 1999; Powers et al., 1998; Seagren et al., 1999). In contrast, other studies have reported a kinetic controlled dissolution of DNAPL (Miller et al., 1998; Nambi and Powers, 2003; Saba and Illangasekare, 2000). In this integrated modeling framework, a flexible numerical formulation for dissolution is defined based on stagnant film theory (Pankow and Cherry, 1996; Equation 7):

$$q = K_l(C_s - c) \tag{eq. 7}$$

$K_l$ is the lumped mass transfer coefficient (m·s$^{-1}$), $C_s$ is the equilibrium concentration of contaminant in water (kg$_{TCE}$·m$^{-3}$).

For kinetically controlled dissolution, numerous laboratory and numerical techniques have been developed to estimate the mass transfer coefficient (Dalla et al., 2002; Gvirtzman et al., 1987; Kim et al., 1999; Ronen et al., 1986). More specifically, the transfer rate has been reported to be dependent on the effective solubilities and diffusivities of the DNAPL components, the DNAPL source zone architecture (Luciano et al., 2018), fluid flow properties (Miller et al., 1990; Nambi and Powers, 2003a; Saba and Illangasekare, 2000), and sand grain properties (Nambi and Powers, 2003a; Powers et al., 1994). Each approach for the estimation of transfer rate has its own complexity and limitation which necessitates field specific estimation (Grant, 2005). Seargen 1994 reported that Local equilibrium Assumption is possible in an aquifer system where significant contact between the flowing groundwater and the NAPL contamination and those in which the interphase mass-transfer rate is relatively large (Seagren et al., 1994). In order to simulate a simpler but plausible case, a Local Equilibrium Assumption (LEA) for dissolution processes is considered in the hypothetical aquifer. For simulating the case of LEA, mass transfer coefficient in Equation 7 can be assigned a significantly higher value that results in the equilibrium governed dissolution. Along with the mas transfer coefficient, The values of all parameters used for the simulating the TCE dissolution is summarized in Table 3.

*Table 3: Dissolution and transport properties of TCE used for Stage 2 simulation in the hypothetical aquifer.*

| Parameters | Value |
|---|---|
| Solbility of TCE | 1.27 Kg·m$^{-3}$ |
| Dispersivity | 0.02 m |
| Dispersion coefficient | 1×10$^{-8}$ m$^2$·s$^{-1}$ |
| Mass transfer rate | 1200 |

### 2.4.1 Initial conditions

Initially the hypothetical aquifer is considered to be free of the dissolved TCE. The source for the dissolved contaminants is numerically prescribed in the region where the model predicts formation of TCE pool and ganglia at the end of stage 1. A dispersivity of 2 cm has been considered for the contaminant plume which lies in the typical range for a real field aquifer. The left and right boundaries of the aquifer have been considered as an open boundary for contaminant flow. The top and bottom boundaries are no flow boundaries.

Unlike in Stage 1, the groundwater flow has been considered to simulate the transport of the TCE. Different hydraulic heads are assigned at the lateral end of the model domain to simulate groundwater flow. While 12 m of end is assigned at the right boundary, 13.25 m of hydraulic head is assigned at the left boundary. The aqueous TCE model is run to demonstrate a long-term contamination of groundwater in decade scale. For this, stage 2 is run for 11 years after which the after which the majority of TCE ganglia gets depleted, while the pool persists in the hypothetical aquifer.

## 2.5 Modelling Stage 3: nZVI injection

### 2.5.1 Initial and boundary condition

In Stage 3, the transport and retention behavior of ZVI in the hypothetical aquifer is simulated. The simulation aims to study the transport and retention behavior of the nZVI including the prediction of its Radius of Influence (ROI). In this stage, the initial and boundary condition of nZVI injection has been set with the intention of a pilot-scale remediation resulting in partial degradation of contaminant. An injection well is placed at the downgradient of the expected TCE source zone in the upper sand layer. A well screen is considered at a depth of 6.6 meters for the injection of CMC-nZVI in the simulated domain (Figure 3). The model has been developed with the aim to simulate the degradation of the contaminants at the upper sand. For simulating the nZVI fate, the model requires several parameters such as injection rate, CMC concentration, ZVI concentration.

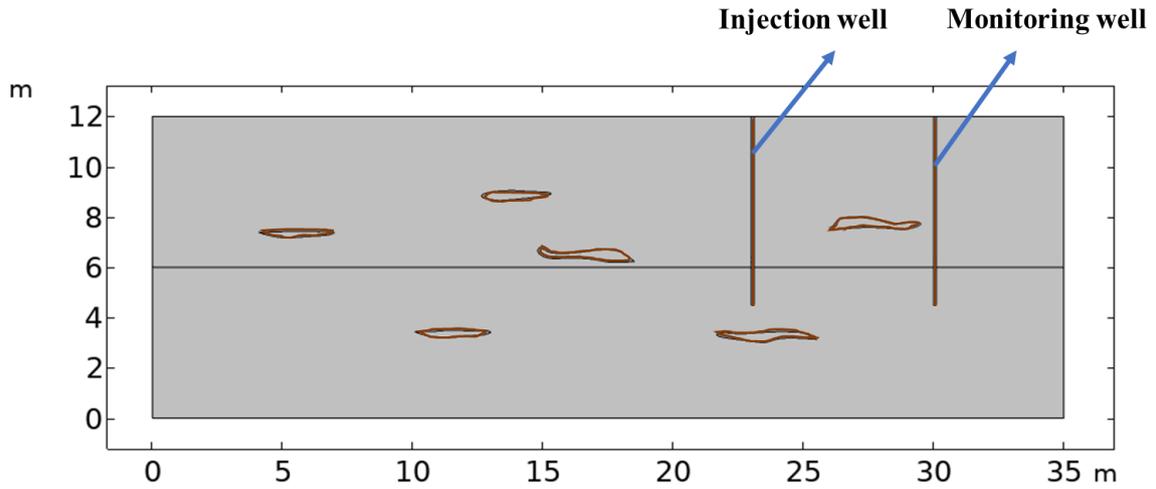

*Figure 3: Model domain with the incorporation of the injection and monitoring well placed 7 meters apart. The 2cm well screen for both injection and monitoring well is assigned at a depth of 6.6 meters from the surface,*

A modified advection dispersion retention equation has been used for the simulation of nano-particles in the groundwater (Equation 6, 7):

$$\frac{\partial(\theta c)}{\partial t} + \rho_b \frac{\partial(s)}{\partial t} = -\nabla(vc) + \nabla(D\nabla c) \qquad \text{(eq. 6)}$$

$$\rho_b \frac{\partial(s)}{\partial t} = \theta k_a c \qquad \text{(eq. 7)}$$

Equation 6 represents the mass conservation of the aqueous nZVI in the hypothetical aquifer, whereas Equation 7 represents the first order deposition of nZVI. The retention of the ZVI is characterized by Colloid Filtration Theory where the deposition of particle is considered only in the favorable sites (Yao *et al.*, 1971). The theory depicts the attachment of nano-particles in the sand grain due to Brownian diffusion, interception, and gravitational sedimentation (Tufenkji and Elimelech, 2004; Yao *et al.*, 1971). This equation does not take into account particle-particle interaction or detachment of the deposited nZVI. Previous researchers have reported that the ZVI can be efficiently simulated with Colloid Filtration Theory (Asad *et al.*, 2021; Krol *et al.*, 2013). The formulation for attachment rate of the particle is represented by Equation 9:

$$k_{att} = \frac{3}{2} \frac{(1-\theta)v\alpha\eta_o}{d_c} \qquad \text{(eq. 9)}$$

where $d_c$ is the collector grain size (diameter of the aquifer grain) (m), $v$ is the pore water velocity (m·s-1), α is the attachment efficiency (-) which characterizes the ratio of nanoparticles deposited on the collector per nanoparticle collision, and $η_o$ is the the the single collector contact efficiency (-), which quantifies the ratio of particles that strike the collector to those that approach the collector. $η_o$ was estimated based on the dimensionless parameters that are associated with the phenomena of diffusion, interception, and gravity (Tufenkji and Elimelech, 2004, Krol et al., 2013). The parameters corresponding to each phenomena were determined based on the mathematical relationship established in a previous study from the literature (Tufenkji and Elimelech, 2004).

For this hypothetical aquifer, the values of this parameter is taken from the previous CMC-nZVI injection case at a field site near San Francisco Bay, as reported in the literature (Krol *et al.*, 2013). The values of the parameters such as nZVI diameter, sand diameter (collector diameter), injection rate, and nZVI concentration are listed in the Table 4.

*Table 4: Values of the parameter used for the simulation of CMC-ZVI´s mobility in the hypothetical aquifer*

| Parameters | Value |
|---|---|
| Concentration of injected nZVI | 0.2 $kg_{nZVi} \cdot m_{water}^{-3}$ |
| Darcy velocity for injection | 88 $m \cdot day^{-1}$ |
| nZVI particle diameter | 140 nm |
| Attachment efficiency | 0.02 |
| Polymer solution viscosity | 0.0027 Pa·s |
| Polymer molar Fraction | $8.04 \times 10^{-7}$ |
| Concentration of injected CMC | 3 $kg_{CMC} \cdot m_{water}^{-3}$ |

The size of the sand grain (collector diameter $d_C$) is expected to vary with the permeability field.. The Kozeny-Carman equation(Tsakiroglou et al., 2018) was used to derive the value of collector´s diameter based on the values of the permeability and porosity of different lithological units within the hypothetical aquifer (Equation 8):

$$d_C = \sqrt{\frac{k_p(1-\theta)^2 180}{\theta^3}} \qquad (eq.\ 8)$$

### 2.5.2 Formulation for clogging of nZVI in porous media

The retention of nZVI in the porous media is expected to alter the hydrodynamic properties of the host mineral grain due to the clogging ef(Tsakiroglou *et al.*, 2018). In order to estimate the change in porosity and permeability, the model incorporates the formulation used in previous studies (Tosco et al., 2009, 2014a, 2014b; Equation 9, 10, 11)

$$\theta_m = \theta_0 - \frac{\rho_b}{\rho_s} s \qquad (eq.\ 9)$$

$$a(s) = a_0 + a_p \gamma \left(\frac{\rho_b}{\rho_p}\right) s \qquad (eq.\ 10)$$

$$k(s) = k_0 \left(\frac{\varepsilon_m}{\varepsilon_0}\right)^3 \left(\frac{a_0}{a}\right)^2 \qquad \text{(eq. 11)}$$

In Eqs. (9)–(10), $\theta_0$ is the initial porosity, $\theta$ is the fraction of initial porosity occupied by deposited particles, $\rho_s$ is the density of the deposited particles, $a_p$, is the specific surface area of the ZVI, $a_0$ is the initial specific surface area of the porous medium, $k_0$ is the initial absolute permeability of the porous medium, and $\gamma$ is the parameter representing the fraction of deposited particles contributing to the increase of the surface area of the matrix. The value of the parameters for simulating the clogging effect of ZVI has been taken from the previous report (Tosco and Sethi, 2010) and has been summarized in Table 5.

*Table 5 parameters used for the simulation of clogging which is taken from a previous research study (Tosco and Sethi 2010)*

| Parameters | Value |
|---|---|
| Specific surface area of silty sand | $4.99 \times 10^3$ (m$^{-1}$) |
| Specific surface area of zvi | $2.34 \times 10^8$ (m$^{-1}$) |
| Density of ZVI | 6100 (kg/m$^3$) |
| Densty of sand | 2600 (kg/m$^3$) |
| Fraction of ZVI altering surface area ($\gamma$) | $1.04 \times 10^{-3}$ |

### 2.5.3 TCE plume disruption due to nZVI injection

The injection fluid is expected to disrupt the natural plume geometry of aqueous TCE during the course of nZVI injection. To study the effect of plume disruption, the transport model of aqueous TCE is run simultaneously along with the simulation for nZVI injection. During the injection period, the evolution of TCE source zone architecture is also simulated to update the secondary source during the injection. Under the integrated modeling framework, the output from stage 2 for aqueous TCE plume geometry and TCE source zone architecture are considered as the initial condition for this simulation.

### 2.5.4 Effect of viscosity change

While CMC helps in the efficient transport of the nZVI, it also alters the viscosity of the fluid medium, thus influencing the flow field (Krol et al., 2013; Li et al., 2015). Therefore, in addition to the simulation of the nZVI´s mobility, the model simulates the effect of CMC (injected along with nZVI) on the viscosity of the groundwater. CMC was modeled as a conservative species as previously done in the literature reports (Asad et al., 2021; Krol et al., 2013). The result of the CMC plume simulation is used to predict the evolution of viscosity in the groundwater. For this, Grunberg and Nissan equation (Equation 12) is used which represents viscosity change as a function of CMC concentration.

$$ln(\mu_{sol}) = x_{CMC} \, ln(\mu_{CMC}) + x_{water} \, ln(\mu_{water}) \qquad \text{(eq. 12)}$$

Where $\mu_{sol}$ is the viscosity of the aqueous phase carrying CMC, $\mu_{CMC}$ is the viscosity of the injected fluid, $x_{CMC}$ is the mole fraction of CMC, and $x_{water}$ is the mole fraction of water. The concentration of injected CMC has been set at 3 kg$_{CMC}$·m$_{water}^{-3}$ to resemble the experimental configuration reported in a previous study (Krol et al., 2013).

The nZVI injection is simulated for 8 hours after which the model demonstrates a sufficient Radius of Influence for the degradation of simulated TCE plume. In Stage 4, a total of five different phenomena are simulated, namely the disruption in the TCE transport, the TCE saturation evolution, the conservative transport of CMC, fate of nZVI, and the clogging effect due to nZVI injection. The significant output of this model demonstrates the modified groundwater flow field due to clogging of nZVI and viscosity change in the presence of CMC, disruption in TCE plume due to injection and importantly, the spatial distribution of retained ZVI with a estimated ROI. After simulating the injection of nZVI injection, the groundwater in the model domain has been set at its natural flow condition.

## 2.6 Modelling Stage 4: Reactive transport simulations

The simulation in Stage 4 focuses on the study, the efficiency, and timeline for the remediation with nZVI. Therefore, a reactive transport model is developed to simulate the degradation of the aqueous TCE plume in presence of nZVI. Simultaneous to the degradation study, the model simulates the surface passivation of nZVI surface due to its reaction with TCE. The coupled simulation aims to provide predictive insight into the fate of TCE and the longevity of nZVI. Under the integrated modeling framework, the input for Stage 4 corresponds to the output of Stage 3 after 8 hours of CMC-nZVI injection.

Several researchers have reported different reactive pathways, kinetic rates and the efficiency of nZVI for the degradation of TCE. While multi-reaction pathway governs the TCE degradation in the real aquifer, acetylene is predominant degradation product where the TCE concentration is relatively higher(Liu *et al.*, 2005). As the model is based on LEA resulting in TCE concentration up to its solubility limit, acetylene is assumed to be the only degradation product. The reaction for TCE degradation in the presence of nZVI is represented by Equation 13 (Prommer *et al.*, 2008):

$$2Fe^0(s) + C_2HCl_3 + H+ \Rightarrow C_2H_2 + 2Fe^{+2} + 3Cl^- \qquad \text{(eq. 13)}$$

It has been widely observed that the TCE degradation follow pseudo first order kinetic reaction under batch condition(Gomes, 2014; Song and Carraway, 2005; Xu *et al.*, 2019). The formulation for the pseudo-first order reaction rate equation for the TCE degradation is represented by equation 2.65:

$$-\frac{dc}{dt} = k_{SA} \propto_s \rho_m c \qquad \text{(eq. 14)}$$

Where $\propto_s$ is the surface area m²·g⁻¹, and $\rho_m$ is the concentration of the iron distributed in the porous media $kg_{nZVI} \cdot m_{water}^{-3}$. The $k_{SA}$ value depends on the nZVI characterization and for this research work, $k_{SA}$ a value of $2.6 \times 10^{-3}$ L·h⁻¹·m⁻² and surface area ($\propto_s$) of 23 m²·g⁻¹ are taken from the literature (Liu *et al.*, 2005). Considering the degradation reaction, the set of governing equations for the mobility of contaminants in the aquifer has been represented by the mass conservation expression (Equation 5.2):

$$\theta \frac{dc}{dt} + v\nabla c - \nabla D(\nabla c) = K_l(C_s - c) - k_{SA} \propto_s \rho_m c \qquad \text{(eq. 15)}$$

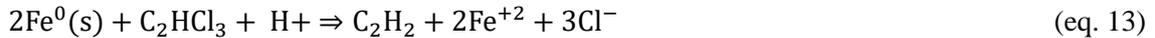

The TCE simulation assumes the accessibility of core shell structure of nZVI with negligible effect of oxidation layer formation due to corrosion. The depletion of active surface of nZVI due to its reaction with TCE is represented by Equation 16. The equation considers the stoichiometric ratio of TCE and iron surface explained in Equation 13.

$$\frac{d\rho_m}{dt} = -0.85 k_{SA} \propto_s \rho_m c \qquad (eq.\ 16)$$

The model is simulated for 6 years to predict the degradation of TCE in the presence of nZVI retained in the lithological units of the hypothetical aquifer.

Table 6 summarizes the coupled simulations implemented in the integrated modeling framework. The TCE phase infiltration (Stage 1) in the model domain has been simulated for 135 days which predicts the TCE source zone architecture. Subsequently in Stage 2, the simulation of dissolution and transport of aqueous TCE is added in the integrated framework. The model is implemented with the objective to study the spatio-temporal evolution of TCE plume for 11 years. In Stage 3, the model for characterizing the fate of nZVI is further added in the integrated modeling framework. In the last stage, contaminant degradation is simulated for the following 2.5 years after the injection of nZVIs. The last stage simultaneously simulates the evolution of ROI of nZVI as well as aqueous and source zone of TCE. As the values are taken from different literature, it is highly likely that the hydrogeological and nZVI properties for a real case deviate from the hypothetical case considered here. However, the purpose of this research is limited to the demonstration of the capabilities of end-to-end modeling framework in providing predictive insights for a comprehensive remediation design. The verification of numerical simulation used for the end-to-end numerical simulation has already been done in a previous research work (Nilabh, 2021)

*Table 6: List of the different combination of simulation used in the different stages.*

|  | TCE source zone geometry | TCE dissolution and Transport | Radius of Influence of nZVI | TCE degradation |
|---|---|---|---|---|
| Stage 1 (135 days) | ✓ | - | - | - |
| Stage 2 (11 years) | ✓ | ✓ | - | - |
| Stage 3 (8 hours) | ✓ | ✓ | ✓ | - |
| Stage 4 (2 years) | ✓ | ✓ | ✓ | ✓ |

# 3 Results

## 3.1 Stage 1: Two-phase Flow of TCE infiltration in groundwater

In stage 1, the simultaneous two-phase flow of TCE and water is simulated to predict (a) the migration pathways of infiltrated TCE, (b) the travel time of TCE before immobilization, and (c) the TCE architecture at the source zone. The Figure 4 shows the model prediction for the TCE migration pathway after 5 days, 15 days, 25 days and 35 days of simulated TCE migration. For a better visualization an enlarged plot focusing on the TCE mass is obtained. Furthermore, the plotting for saturation profile is done only for the region where the TCE saturation is non-zero.

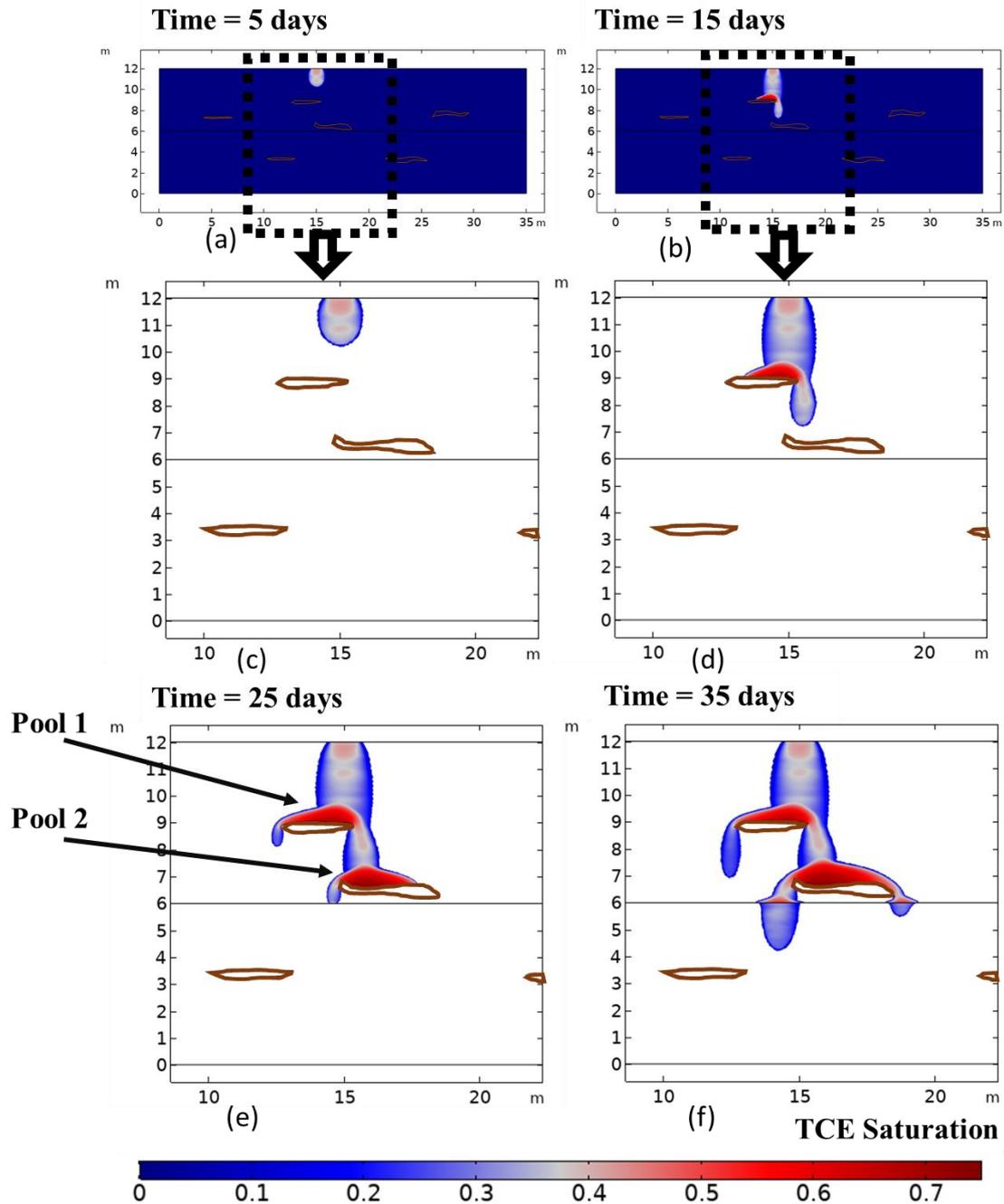

*Figure 4: Model prediction for TCE saturation during the Stage 1. (a) Predicted TCE plume after 5 days of infiltration. (b) after 15 days of infiltration, and their zoomed-in image focusing on TCE phase migration (c), (d) Likewise, TCE saturation profile after (e) 25 days and (f) 35 days of infiltration*

At time (t) = 5 days, the predicted TCE plume moves predominantly in the vertical direction with a negligible lateral movement beyond the extent of the infiltration area (Figure 4a). The model further estimates a varied range of TCE saturation with a maximum value of 0.5. The central saturation of the TCE phase tends to have a non-uniform value at different depths. The non-uniform TCE saturation highlights the role of microscopic heterogeneity resulting in the formation of two different patches of relatively higher TCE saturation. After 15 days of the TCE infiltration, the model predicts the accumulation of TCE at the top of the clay lens with a maximum saturation of 0.7 (Figure 4b). As the simulated phase fails to infiltrate into the clay lenses, the result indicates that the entry pressure of upper clay is higher than the capillary pressure at the overlying sand.

After 25 days, the model predicts the TCE vertical motion as the two branches of TCE plume are formed at the two edges of the upper clay lens (Figure 4c). The vertical migration of both the branches of TCE leads to the formation of a second TCE pool over the encountered clay lens. The model estimates a maximum saturation of 0.69 for Pool 1 and 0.75 for Pool 2. Furthermore, the model prediction illustrates the onset of vertical migration of TCE phase from Pool 2. However, the vertical migration of this TCE mass is hindered at the interface of the upper and lower silty layer. The result again highlights the role of entry pressure of the lithology governing the flow of TCE. Eventually, after time (t) = 35 days, the model prediction shows the vertical downward motion of the TCE plume from the four branches each formed at the edges of the two TCE pools (Figure 4d). Similar to the sand-clay interface, the accumulation of simulated TCE mass at the silty sand interfaces is observed. However, unlike the sand-clay interface the model predicts the infiltration of the TCE mass in the low-permeable silty sand. Thus, the results indicate that the capillary pressure of the simulated TCE phase front overcomes the entry pressure of the lower silty sand.

After 35 days, the further downward movement of TCE phase is driven by gravity and capillary forces in the absence of the external TCE source at the top of the hypothetical aquifer. Figure 5 illustrates the distribution of TCE phase saturation for the subsequent days. After 40 days, the model predicts a gradual decrease in the TCE saturation from the top of the hypothetical aquifer (Figure 5). Three TCE phase fronts, PF1, PF2 and PF3 are developed, with an infiltration range up to 5.6 meters, 8.7 meters, and 7.2 meter respectively when measured from the top of the domain. Subsequently, after 60 days, all of these phase fronts are estimated to have infiltrated the lower sand (Figure 5b), with PF1 extending up to the bedrock of the aquifer domain. The saturation of the TCE pools decreases progressively with the maximum estimated saturation of 0.53 at Pool 1 and 0.68 at Pool 2. Compared to the TCE saturation profile after the infiltration for 35 days (when the external source at the top of the model domain is stopped), the decrease in the maximum saturation at pool 1 and pool 2 amounts to 25 % and 9% respectively. The decrease in the TCE saturation at Pool 1 and 2 indicates the downward remobilization of TCE pool in absence of external source. Additionally, the decrease is relatively higher for pool 1 compared to pool 2 suggesting varied degree of remobilization of TCE phase. The result thus indicates the role of clay lens morphology in the accumulation and remobilization of TCE phase. In contrast, due to the downward migration of the phase fronts, three TCE pools with moderate saturation of 0.46 is formed at the interface of the upper and lower silty sand layers. After 85 days, the model estimates the formation of a TCE pool at the bedrock due to the accumulation of PF2 front (Figure 5c). Furthermore, the simulated PF3 plume completes its vertical migration resulting in the onset of TCE pool formation at the bedrock. The prediction indicates that the saturation of the TCE Pool 1 and 2 decreases further and ranges up to 0.47 and 0.66 respectively.

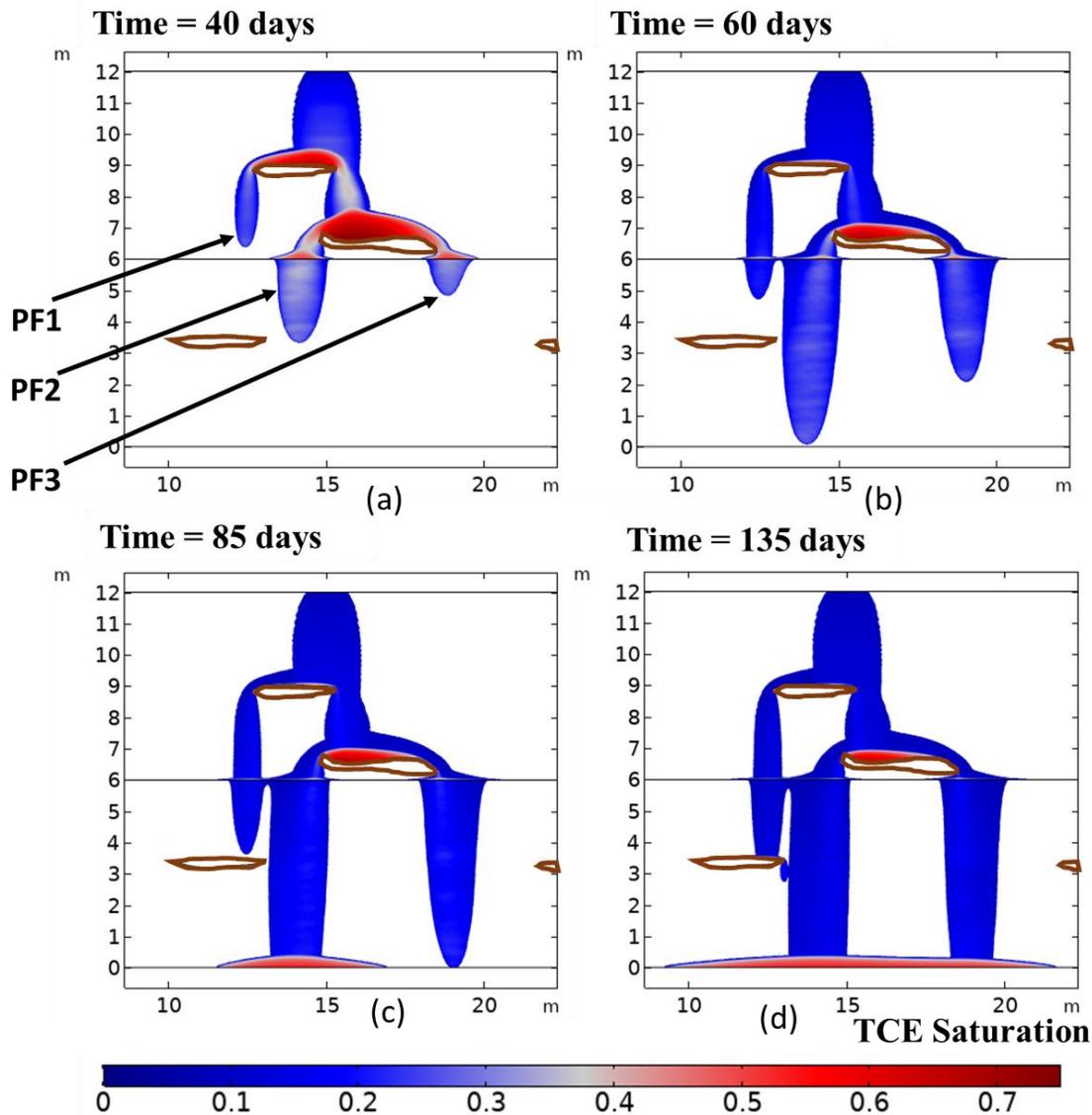

*Figure 5: Model prediction for TCE migration after the removal of external source zone from the top of model domain. (a) simulated TCE phase after 40 days; (b) Prediction of the arrival of TCE plume (PF2) at the bedrock after 60 days of infiltration; (c) TCE saturation profile along the migration path and accumulation at the impervious layer after 85 days. (d) near-static TCE source zone geometry after 135 days*

After 135 days, the model predicts a near-static TCE phase mass comprised with a relatively large TCE pool at the bottom of the aquifer. The model predicts an integration of TCE pool formed by both PF2 and PF3 plume resulting in 13 m horizontal extent of TCE at the bedrock (Figure 5d). The saturation of TCE pool over the clay lenses ranges up to 0.43-0.66 while the saturation of the pool at the bedrock ranges up to 0.53. In contrast, the ganglia saturation is predicted to be in the range of 0.01-0.15

The static TCE source zone consists of ganglia along the infiltration pathway and pool at bedrock and lithological interfaces within the simulated domain. A total of 32.8 % of simulated TCE mass remains in the upper layer, while the remaining 68.2 % remains in lower layer. The entire TCE ganglia constitutes 70 % of the total TCE mass in the aquifer, while the remaining 30 % of TCE are in the form of pool. As the TCE plume becomes nearly static, it is considered as an immobilized TCE source zone. Therefore, the TCE saturation distribution at time (t) = 135 days (Figure 5d) is considered to be the initial TCE source architecture for the Stage 2 simulations under the integrated modeling framework.

## 3.2 Stage 2: Dissolution and transport of TCE

The onset of stage 2 begins with the prediction of groundwater flow due to the pressure head differences assigned at the lateral boundaries. The Figure 6 shows the streamlines of groundwater flow in the simulated domain. The color map in the streamline depicts the estimated groundwater flow speed. The result indicates that the simulated groundwater flow rate is non-uniform and ranges from 1.4 to 5.5 cm·day$^{-1}$. The non-uniform groundwater flow rate predicted by the model demonstrates the role of microscopic heterogeneity governing the flow dynamics. The flow rate is relatively higher in the upper silty sand layer compared to the lower silty sand and can be attributed to their difference in permeability. The flow rate in the lower silty sand varies from 1.4 to 3.8 cm·day$^{-1}$. The arrow in Figure 6 implies that groundwater movement is predominantly horizontal. However, close to the clay lenses, the flow has local vertical component due to the obstruction from low permeability zone. The groundwater streamline near to the bedrock illustrates a slightly elevated flow field where the TCE pool formation has been predicted in Stage 1. The result thus highlights the influence of TCE saturation in lowering the permeability of the porous media and thus groundwater flow rate.

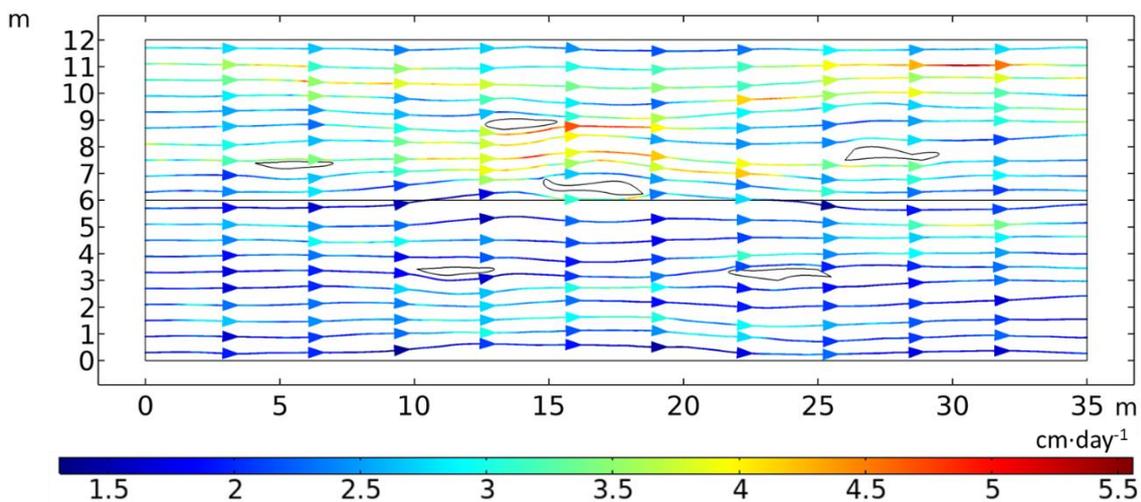

*Figure 6: Groundwater velocity field for the Stage 2. The Darcy velocity at the upper silty sand layer ranges up to 5.5 cm·day$^{-1}$ while the Darcy velocity in lower sand ranges up to 3.5 cm·day$^{-1}$*

The model further simulates 11 years of TCE phase dissolution and its transport in the aquifer. Consequently, the model results demonstrate the spatio-temporal evolution of aqueous TCE plume and the depletion of TCE phase acting as a source zone (Figure 7). The result for the aqueous TCE transport shows that within 0.1 year, the TCE aqueous concentration reached the solubility limit of 1.27 kg$_{TCE}$·m$^{-3}$$_{water}$ (Figure 7c). The geometry of the simulated aqueous TCE plume resembles with the source zone (Figure 7a) indicating the TCE dissolution from the entire source zone. The maximum spatial extent of aqueous TCE transport is predicted to be 4 m from the source zone. The dissolution during the first 0.1 year is minimal with 97% of the simulated TCE source zone remain undissolved (Figure 7a). The smooth zig-zag shaped front of the simulated aqueous TCE plume indicates that the advection is non-uniform due to the heterogeneous groundwater field. The plume geometry front, nevertheless, is predicted to be smoothened by the hydrodynamic dispersion of aqueous TCE in the aquifer. The Figure 7d shows the simulation result for the downstream transport of aqueous TCE plume after 0.4 year of dissolution. A relatively more heterogenous and finger-shaped front of aqueous TCE plume is predicted. The maximum spatial extent of the simulated TCE plume reaches the boundary of the hypothetical aquifer with a concentration of 0.11 $kg_{TCE} \cdot m^{-3}_{water}$. The model result further depicts that towards the downstream direction, TCE plume has higher spatial extent in the lower silty sand despite its relatively lower permeability. The result can be attributed to the TCE source zone geometry which is distributed in a higher amount in the lower layer compared to upper layer. After 1 year of the TCE dissolution, the model predicts a uniform TCE concentration with solubility limit concentration in the downstream direction (Figure 7e). The model prediction for the TCE source zone shows that the depletion of TCE is predicted from the upgradient direction (Figure 7b). The TCE source zone direction in the downgradient direction is predicted to remain undissolved. The model estimated a total of 87.5 % of TCE mass source that remains undissolved after 1 year. It can be inferred that the geometry of the contaminant plume front is governed by the hydraulic properties of the aquifer and TCE source zone architecture.

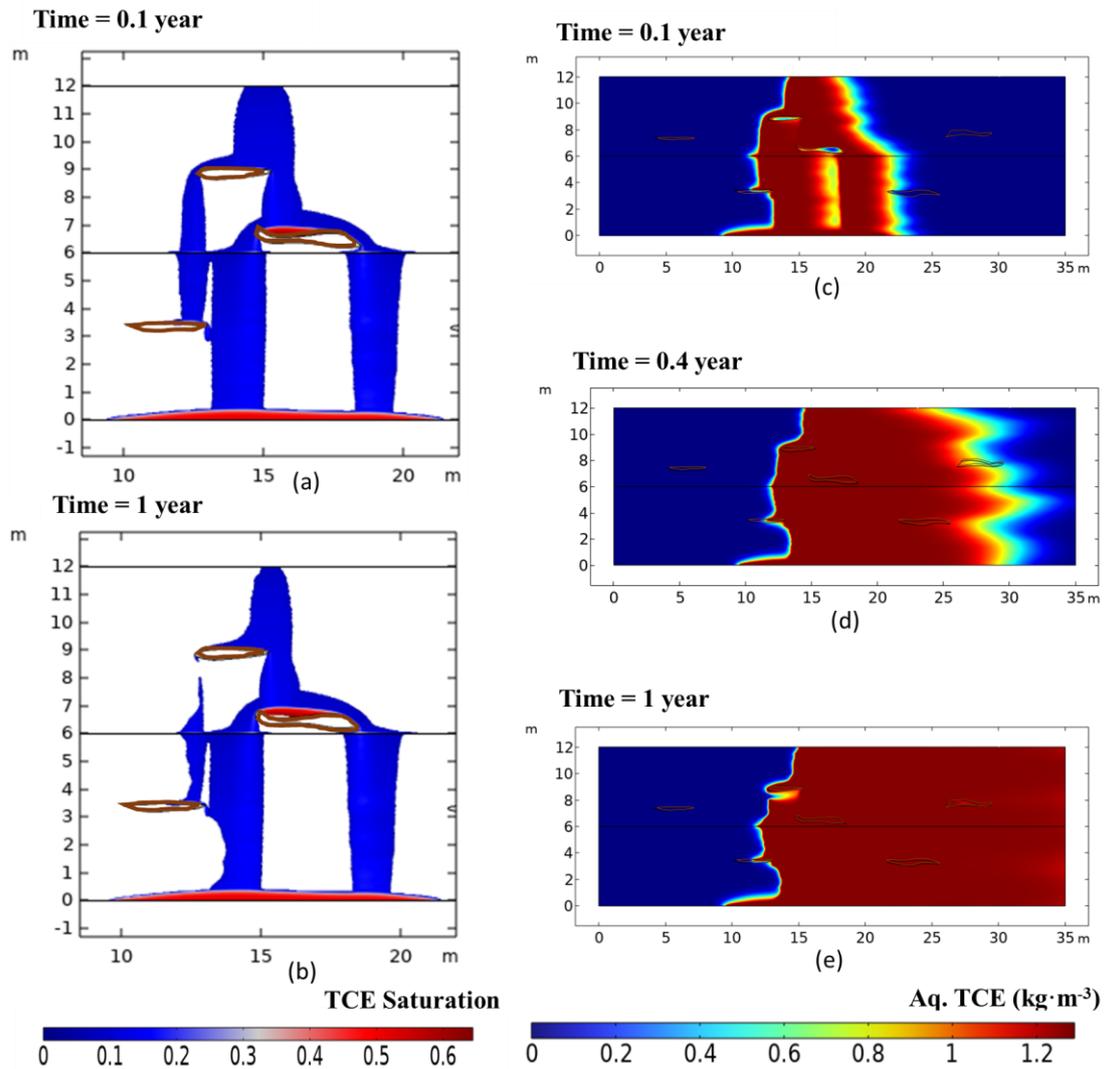

*Figure 7: Model prediction of the evolution of TCE saturation after (a) 0.1 year (b)1 year of dissolution. Simultaneous prediction of aqueous TCE plume geometry after (c) 0.1 year, (d) 0.4 year and (e) 1 year of Stage 2 simulation run*

The model prediction of the TCE transport at longer times (t) = 3 years, (t) = 6 years is obtained. (Figure 8). The result shows that after 3 years of dissolution, 61% of the initial source zone is predicted to remain undissolved (Figure 8a). The simulated TCE ganglia in the upper layer of the domain undergoes substantial depletion with the remaining undissolved mass representing only 6.2% of the source zone after 3 years. In contrast, a minor dissolution of TCE ganglia is predicted at the lower silty sand layer as the undissolved TCE ganglia constitutes 49% of the entire TCE mass after 3 years of dissolution. The varied degree of TCE ganglia dissolution for the upper and lower silty sand can be attributed to the different groundwater flow rate and initial source zone geometry. The result further illustrates that the TCE ganglia of the lower silty sand in the downstream direction remains intact while the dissolution is predicted for the ganglia in the upgradient direction. The model estimates a moderate depletion of TCE pools, with the entire pool mass representing 44% of the total source zone. The TCE source zone after 3 years results in a unique aqueous TCE plume geometry in the simulated domain (Figure 8b). The result shows that in the downgradient direction, the aqueous TCE remains at solubility limit uniformly in the lower sand. The simulated aqueous TCE plume in the upper sand originates primarily from the TCE sources at the top of clay lenses. After 6 years the model predicts that 38% of the initial TCE source zone remain undissolved (Figure 8c). The TCE ganglia in the upper sand as well as the TCE pool of the upper most clay lens (Pool 1) is depleted completely. An enlarged version of Figure 8c is plotted for better visualization where only the non-zero TCE concentration and phase saturation of aqueous and undissolved TCE respectively (Figure 8e). In the lower sand, a significant mass of downgradient TCE ganglia remains undissolved and constitutes 28.2% of the entire TCE mass. In contrast, only a minor mass of the upgradient ganglia remains undissolved and represents only 4.2% of the entire TCE mass. The unique TCE source zone geometry can be attributed to the relatively slower groundwater flow in the lower sand along with a higher mass of TCE source. The model predicts the TCE dissolution from the source zone resulting in a unique aqueous TCE plume geometry (Figure 8 d and f). The simulated aqueous TCE plume in the upper sand originates from Pool 2 and is restricted to a depth below 3.5 m from the surface. Furthermore, the model predicts the uniform concentration of downgradient aqueous TCE in the lower sand originating from the ganglia and the TCE pool. The result from the integrated modeling framework highlights that TCE plume geometry in the upgradient direction evolves with the evolution of TCE source zone.

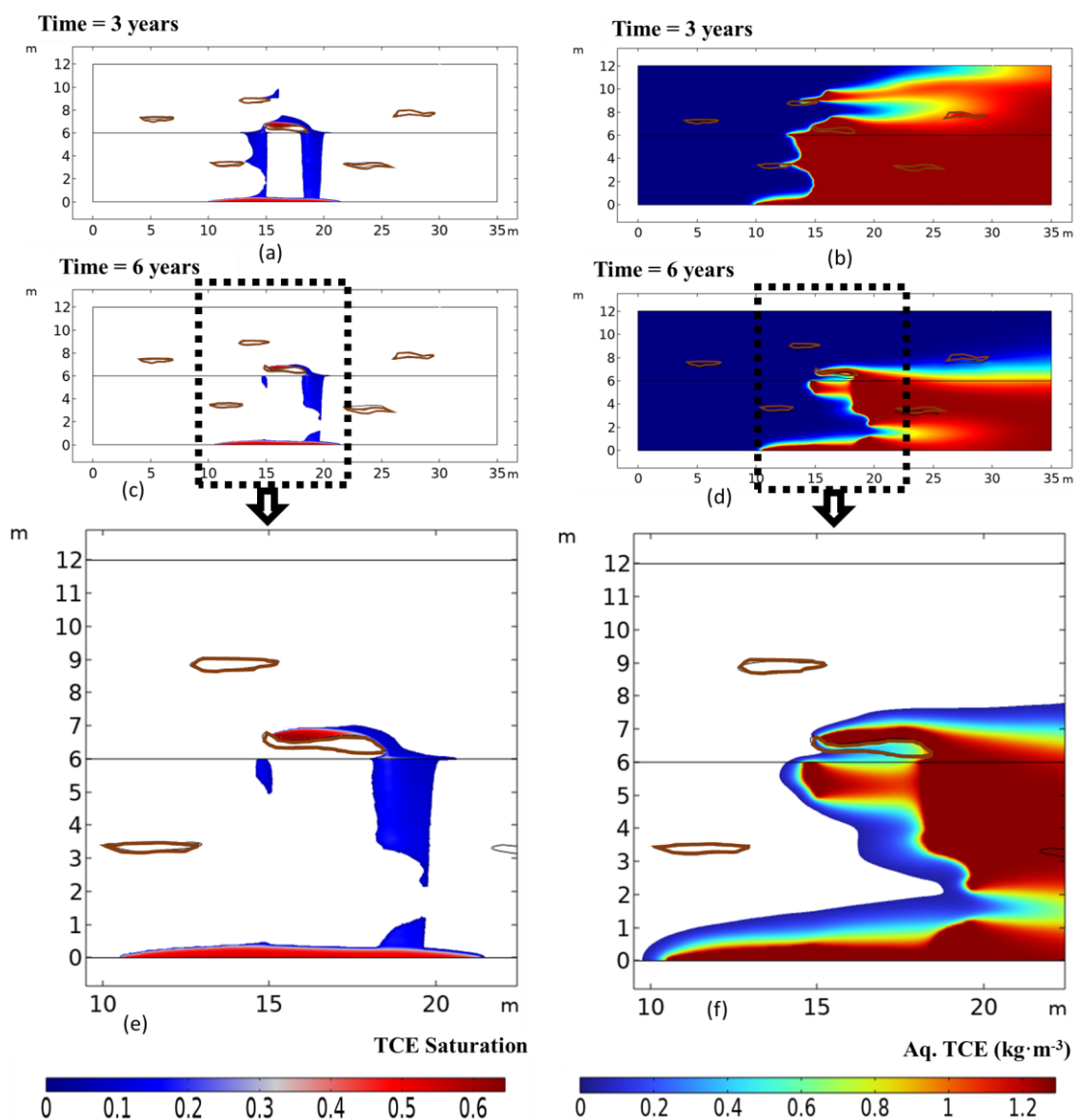

*Figure 8: Model prediction of the evolution of TCE saturation after (a) 3 years (c) 6 years (e) enlarged TCE saturation profile after 6 years where the predicted TCE phase is focused. Simultaneous prediction of aqueous TCE plume geometry after (b)3 years (d) 6 years and (f) enlarged aq. TCE profile after 6 years which focuses on the simulated aq. TCE plume.*

After a period of 11 years of TCE dissolution, the model predicts that only 21.4% of the initial TCE mass remains undissolved (Figure 9b). In addition, the model result indicates that only a TCE pool at upper clay (Pool 2) and bedrock remains at time (t) = 11 years. The enlarged image of the clay lens in the upper sand shows the contamination in the form of TCE pool on the top of the clay lens (Figure 9a). In addition, a small remnant of TCE is present at the upper-lower silty sand interface, resulting in a minor source of aqueous TCE from the interface. The source zone results in the formation of aqueous TCE plume in both upper and lower silty sand(Figure 9d). The enlarger image of the origin of TCE plume highlights the role of shallow level TCE mass in longevity of the contamination in the model domain (Figure 9b). Due to the vertical dispersion of aqueous TCE, the model predicts the concentration of aqueous TCE reaching at the aquifer boundary is lesser than the solubility limit of TCE. The concentration of TCE estimated at the boundary ranges up to 0.47 $kg_{TCE} \cdot m^{-3}_{water}$ and 0.56 $kg_{TCE} \cdot m^{-3}_{water}$ for the shallow level and deeper level TCE plume respectively. The remediation scheme with nZVI is simulated to degrade the simulated plume of shallower level. Under the integrated modeling framework, the starting point of simulation for Stage 3 considers the final TCE phase and aqueous TCE plume geometry obtained in stage 2 simulation (Figure 8 e and f).

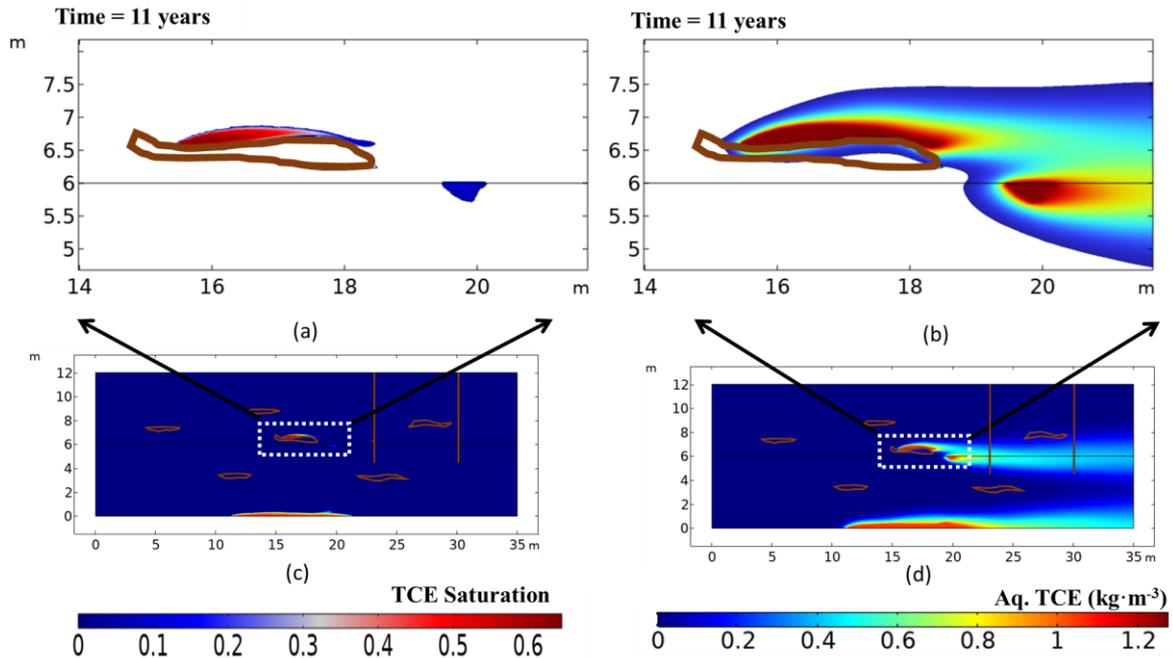

*Figure 9: Initial condition for the modelling of the Stage 3, which corresponds to 11 years since the onset of dissolution of TCE in Stage 2. (a) enlarged image for the TCE pool and ganglia at the upper sand of the simulated domain (b) enlarged image for the shallower aqueous TCE plume at the upper sand of the simulated domain (c) TCE saturation profile (d) aqueous TCE plume geometry in the simulated domain*

## 3.3 Stage 3: Injection of nZVI

### 3.3.1 Flow field generation

The injection of CMC-ZVI within the hypothetical aquifer leads to the generation of a specific flow field (Figure 10). In this flow field, the water originates from the well screen and flows in all directions. The Darcy flow is estimated to be highest at the well screen with a value of 88 $m \cdot day^{-1}$. The flow speed rapidly decreases away from the well screen. While the groundwater flow direction remains the same downgradient with respect to the well, the water flow direction in the upgradient region reverses. The simulation result indicates that the well injection dominates the groundwater flow in the entire hypothetical aquifer. The simulated groundwater flow during the CMC-nZVI injection is used for the simulation of the transport studies of nZVI, CMC and aqueous TCE.

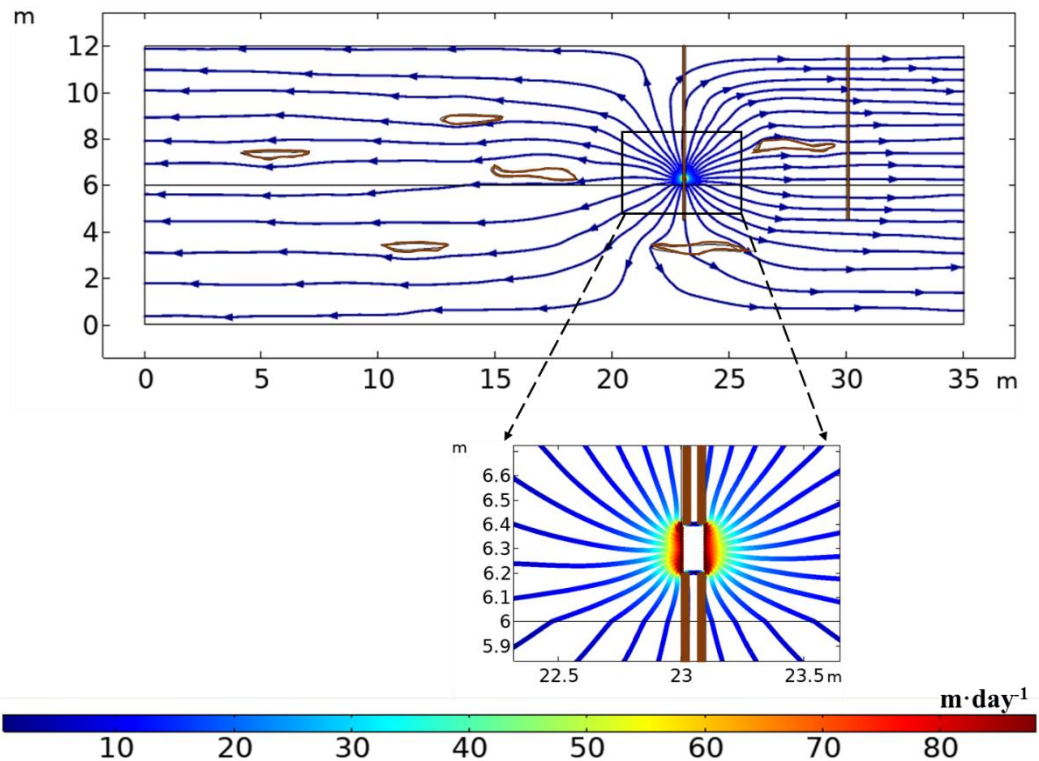

*Figure 10: Model prediction for flow field generated within the simulated aquifer The result highlights that the magnitude of Darcy velocity is highest at the vicinity of the well screen and decreases rapidly with the distance from the screen. The arrows in the plot represent the flow direction of the groundwater.*

### 3.3.2 Impact of the injection fluid on the contaminant plume flow and the effect of viscosity

The onset of injection of CMC-ZVI in the hypothetical aquifer is expected to impact the simulated TCE plume. Therefore, the disruption in the aqueous TCE plume due to the CMC-nZVI injection is simulated under the integrated modeling framework. Figure 11 shows the model prediction for the distorted TCE plume geometry impacted by the injection fluid with CMC-ZVI suspension. After 1 hour of CMC-nZVI injection, the model predicts the development of ring like geometry of TCE plume that envelops the injected water(Figure 11a). The radius of this ring-like geometry ranges upto 0.9 meter. Subsequently, after 4 hours of CMC-ZVI injection, the model predicts a further distortion of TCE plume (Figure 11b). The maximum distance of the predicted TCE plume from the well screen is estimated to be 1.8 meters. At the end of 8 hours, the model predicts the disruption of aqueous TCE plume with an average distance of 4.5 meters from the well screen (Figure 11c). While the aqueous TCE plume is disrupted after the 8 hours injection of CMC-ZVI, the model predicts a steady dissolution of the TCE phase from the source zone as illustrated in Figure 11d. The result indicates that the shape of the aqueous TCE plume is primarily governed by the new groundwater flow field, which depends on the injection set-up and hydrogeological property of the hypothetical aquifer.

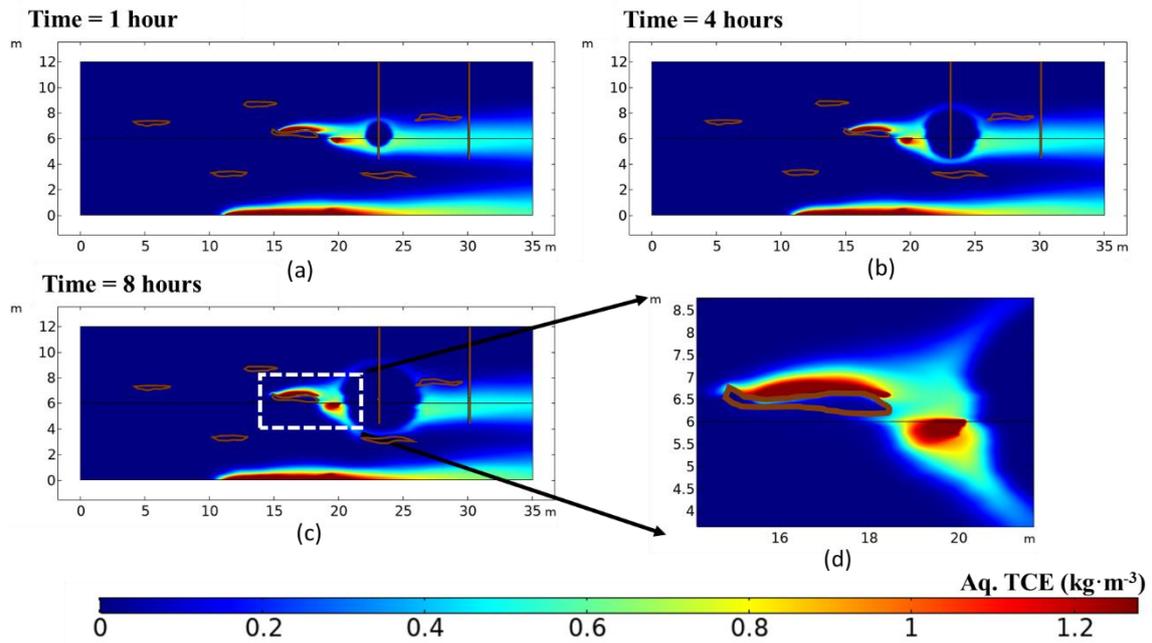

*Figure 11: Model prediction of the impact of the injection of water on the TCE plume (Stage 3), (a) after 1 hour, (b) after 4 hours (c) after 8 hours, and (d) after 8 hours; enlarged around the shallower level TCE source zone*

### 3.3.3 CMC transport and viscosity alteration

The change in viscosity due to the injection of CMC-nZVI solution in the hypothetical aquifer is studied. After 8 hours, the extent of CMC and viscosity change reaches up to 4.5 meters from the well (Figure 12c, d). For the CMC plume with the concentration of 3 kg·m$^{-3}$, the model predicts a corresponding viscosity of 0.0027 Pa·s (Figure 12). The plume geometry of CMC and the viscosity distribution are the same due to their log-linear correlation given by equation 12. The spatial extent of the CMC concentration and viscosity are in tandem with the aqueous TCE disruption studied in the previous section (Section 3.3.2). The figure shows that the geometry of the predicted CMC plume is asymmetrical and thus indicates the role of heterogeneity governing the fate of CMC concentration and subsequently the viscosity change in the hypothetical aquifer. The enlarged image shows highlights that the CMC transport has relatively higher spatial extent for upper silty sand layer compared to lower sand and can be attributed to the difference in their permeability (Figure 12).

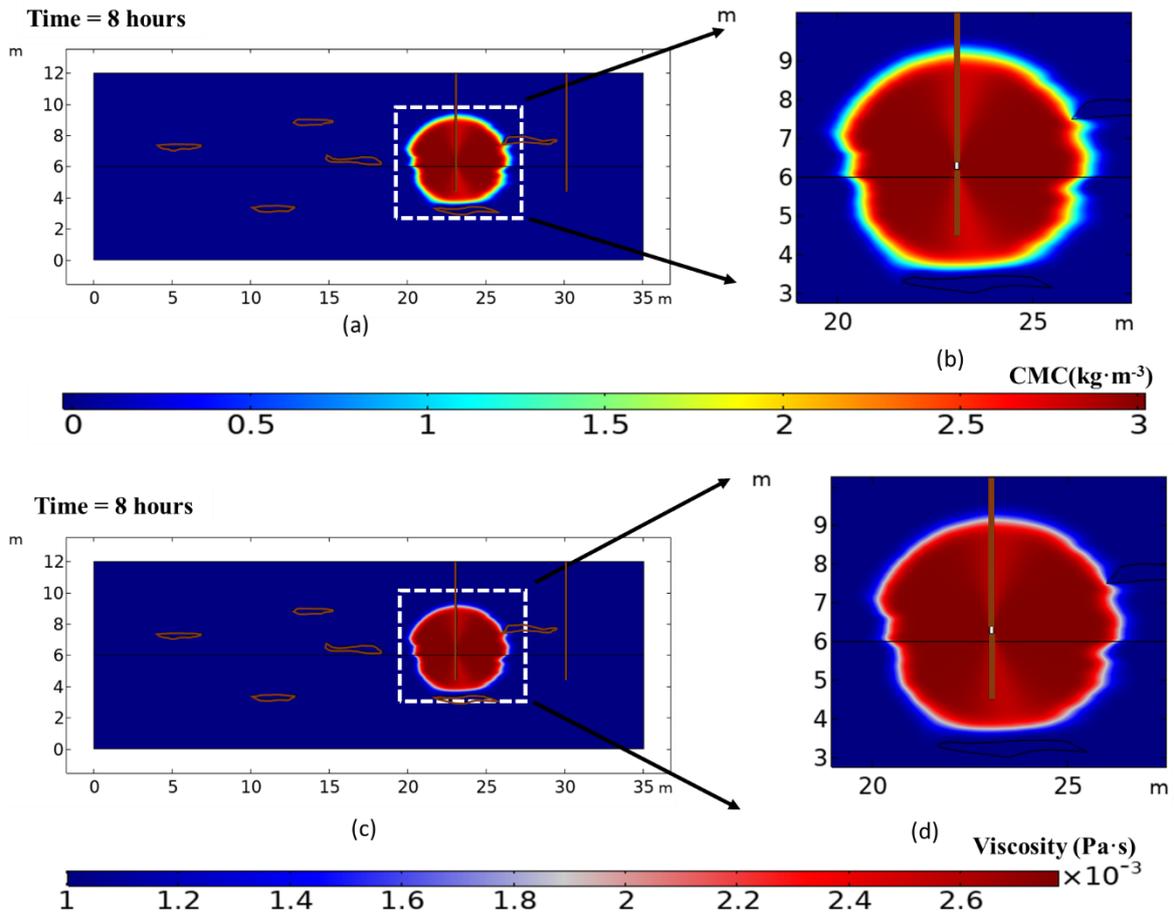

*Figure 12: Model prediction for the CMC concentration due to the fluid injection for (a) 8 hours, (b) 8 hours; enlarged around the injection well screen. Similarly, model prediction of the viscosity within the simulated aquifer after (b) 8 hours of injection. and (d) 8 hours of injection; enlarged around the injection well screen. The profile of viscosity in the simulated aquifer is same as the CMC concentration profile, due to the log-linear relationship between viscosity and CMC concentration.*

### 3.3.4  Injection of nZVI in the hypothetical aquifer

The retained nZVI in the porous media is simulated to study the spatial distribution of the remediating particles within the hypothetical aquifer. The model predicts the dynamic distribution of retained nZVI in the upper and lower silty sand of the hypothetical aquifer (Figure 14). An asymmetrical deposition around the well screen is predicted after 4 hours of injection, with deposition in all directions with respect to the well screen. The simulation result shows the ROI of the retained nZVI is 0.9 meters from the well screen with the highest concentration of 1.53 kg·m$^{-3}$ at the vicinity of the well screen. The enlarged image highlights a contrast in concentration at the interface of upper and lower silty sand (Figure 14). The relatively lower concentration of retained nZVI in the lower silty sand can be attributed to a lower mobility and thus lower retention of nZVI in the relatively less permeable lithological unit. Subsequently, Figure 14c shows the model prediction for the retained nZVI geometry at the end of 8 hours of nZVI injection (Figure 14d). The highest nZVI concentration of 3.14 kg·m$^{-3}$ is estimated at the vicinity of the well screen. The retained ZVI geometry extends up to 1.1 meters for upper sand and 0.9 meters for the lower sand. The higher vertical spatial extent of nZVI in the upper sand layer is favorable to counter the contaminant from the relatively large TCE source in the hypothetical aquifer. Similarly, the smaller spatial extent of deposited particle in the lower sand layer of the model domain is optimized for smaller TCE ganglia. The distance between the predicted TCE source and the estimated rim of the nZVI barrier is 0.9 meters.

Under the integrated modeling framework, the model results after 8 hours of CMC-nZVI injection is considered as the input for stage 4 simulation.

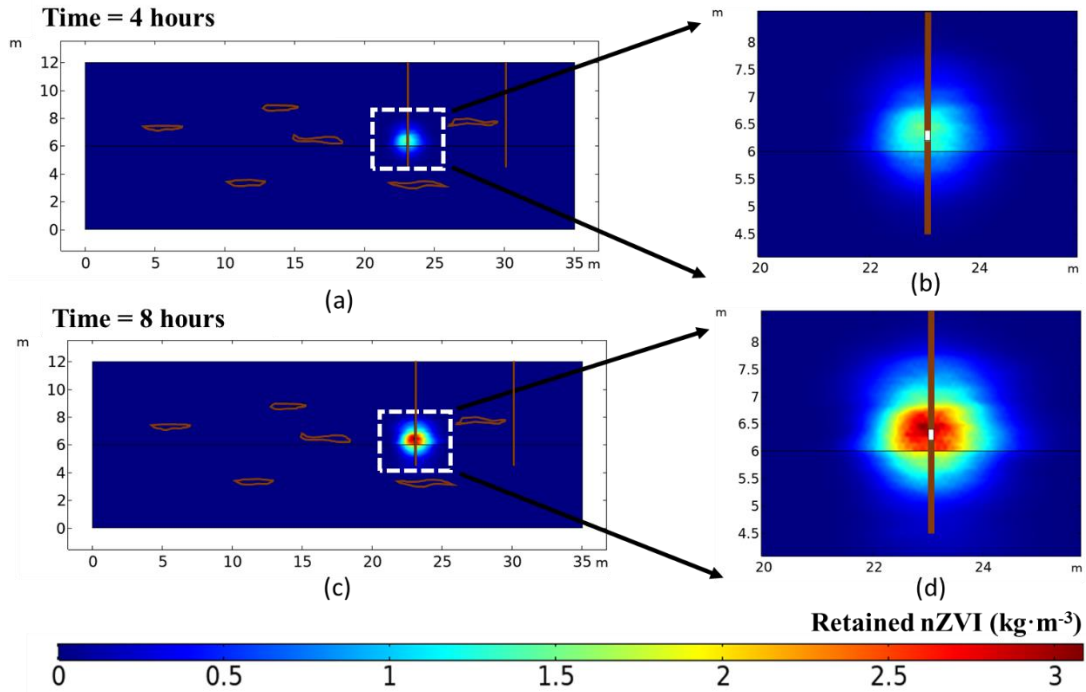

*Figure 13: Model prediction of the retained nZVIs (kg·m⁻³) in the simulated aquifer (a) after 4 hours of injection (b) after 4 hours of injection; enlarged around the injection well screen. aquifer (c) after 8 hours of injection (d) after 8 hours of injection; enlarged around the injection well screen.*

### 3.3.5 Clogging of nZVI in the pores of hypothetical aquifer

The integrated model simulates the clogging effect of nZVI that results in the alteration of porosity and permeability. Figure 14a shows the model prediction for the percentage change in the permeability due to the injection of CMC-nZVI solution in the hypothetical aquifer. The enlarged image for the permeability change shows that a maximum of 4.5% decline in permeability is predicted in the vicinity of the well screen (Figure 14b). The permeability change is higher for the upper sand of the hypothetical aquifer compared to the lower sand. Overall, the shape of the altered permeable area resembles with the spatial distribution of deposited nZVI (Figure 13d). The result illustrates that the alteration of permeability field is primarily driven by the concentration of retained nZVI particles. Figure 14c shows the percentage decrease in porosity after 8 hours of injection. The alteration of the porosity follows the same profile as the permeability with the highest change of 0.12% near the vicinity of the well(Figure 14d).

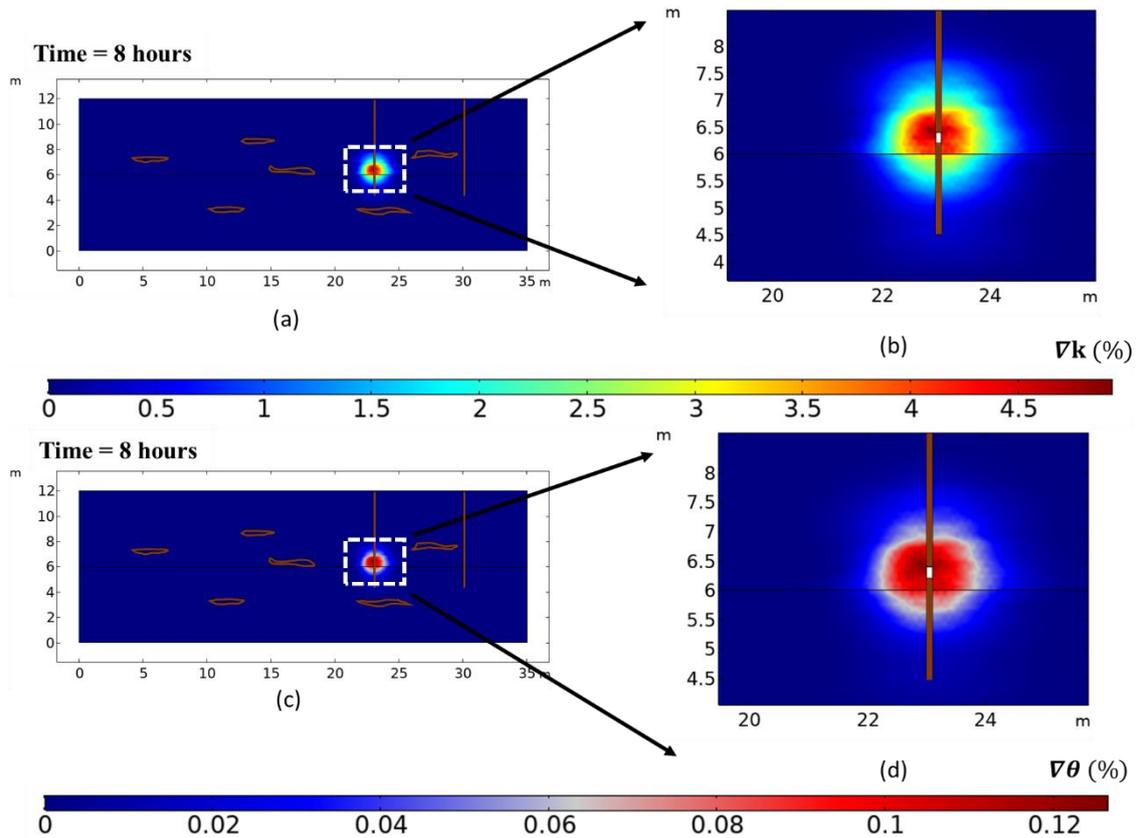

*Figure 14: (a) Model prediction for the % decrease in permeability due to the nZVI injection for the entire 8 hour and, (b) its enlarged version to focus on the well screen. Similarly, (c) Model prediction for the % decrease in porosity due to the nZVI injection for the entire 8 hours and, (b) its enlarged version*

### 3.3.6 Stage 4: TCE degradation in the presence of retained nZVI

The reaction between TCE and nZVI is simulated to study the efficiency of the nZVI in the degradation of TCE. Figure 15a shows the concentration of both aqueous TCE and ZVI in the hypothetical aquifer, after 5 days of their interaction. For a better visualization of the result, only the non-zero values of aqueous TCE and nZVI concentration are shown. The result illustrates that the simulated TCE plume tends to return to the transport driven by natural groundwater flow. The nZVI is visualized with black and grayish gradient with highest concentration represented by grey color at the center. To study the dynamics between the TCE and nZVI, an enlarged image focusing on the interface of nZVI and TCE is plotted (Figure 15d). The figure shows the onset of TCE degradation at the upgradient direction with respect to the nZVI mass. Downgradient to the ZVI mass, the integrated model predicts negligible concentration of TCE and thus no depletion is predicted from that side. After 2 months of TCE-nZVI interaction, the model predicts a formation of pure water zone, as indicated by the white space in the downgradient direction with respect to nZVI mass. The result implies that the nZVI successfully stops the progress of dissolved TCE plume in the downgradient direction, resulting in clean groundwater in the simulated domain. In addition, the model considers the reacted nZVI with dissolved TCE to be no longer active for further degradation, so the barrier is progressively exhausted (Figure 15b). The model estimates the depletion of active nZVI surface in the upgradient direction with 81% of the unreacted nZVI still remaining in the hypothetical aquifer. After 6 months of TCE-nZVI interaction, the model predicts the unreacted nZVI to be 44.8% as the TCE plume transgresses in the nZVI mass (Figure 15 c). The model predicts a zone of length 8 meters of pure water in the downgradient direction. Beyond this zone, the TCE is present at a relatively lower concentration of 0.01 $kg_{TCE} \cdot m^{-3}$, which is the remnant of the TCE plume from stage 3. An enlarged image shows the interface of nZVI-TCE interaction after 6 months (Figure 15 e). The result indicates that the interface is not symmetrical, and a complex TCE-nZVI interface is predicted. Thus, the result indicates that the macroscale heterogeneity of the silty sand results in non-uniform degradation of TCE. The integrated model results emphasize that the progression of TCE plume is governed by the groundwater flow, the retained ZVI concentration and reactivity.

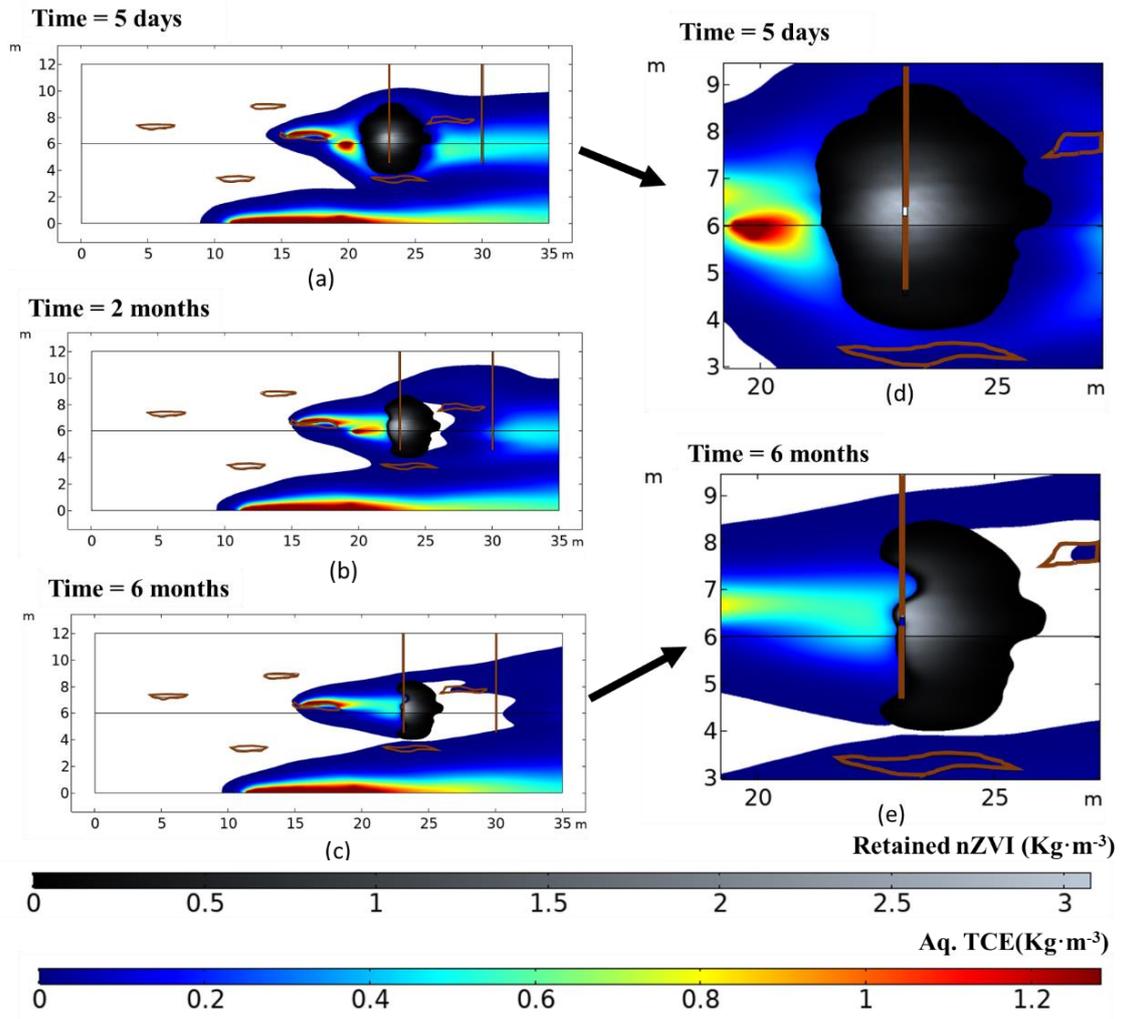

*Figure 15: Model prediction for aqueous TCE plume ($kg_{TCE} \cdot m^{-3}$) water and unreacted nZVI plume ($kg_{nZVI} \cdot m^{-3}$) after (a) 5 days (b) 2 months and (c) 6 months since the injection of nZVI. A rainbow color plot is chosen to represent the areas with non-zero values of aqueous TCE while a grayscale color representation is chosen to represent the nZVI mass. Enlarged version of profile (c) 5 days (d) 6 months since the injection of nZVI, to study the model result at the interface of TCE plume and nZVI mass*

The predicted mobility of nZVI and the contaminant plume in the following years is shown in Figure 16. After 1 year, the model predicts the revival of TCE plume in the hypothetical aquifer (Figure 16a). This transgression through the nZVI mass indicates the formation of favorable channel along the nZVI mass owing to the microscopic heterogeneity. The enlarged version of Figure 16a highlights that the model predicts the bifurcation of active nZVI mass into two different zones (Figure 16d). The total nZVI depletes up to 14.8 % with only 0.02 $kg_{nZVI} \cdot m^{-3}$ as the highest concentration within the unreacted nZVI surface. After 1.5 years, the model predicts the 4.6% of the original nZVI mass to be unreacted (Figure 16b) with a near-complete recovery of TCE plume geometry. The available nZVI mass is localized at the lower silty sand unit. The undepleted nZVI mass in the lithological unit can be attributed to the relatively lower TCE concentration at the rim of TCE plume resulting in slow reaction kinetics. After 2.5 years, the model predicts the complete recovery of plume geometry which resembles with the plume geometry before the injection (Figure 16c). The enlarged image of the figure highlights that an insignificant but 1.1% of the original nZVI mass is still predicted to be undeleted (Figure 16e). However, due to the low concentration of the nZVI mass coupled with the low concentration of TCE results in a steady-state TCE plume geometry.

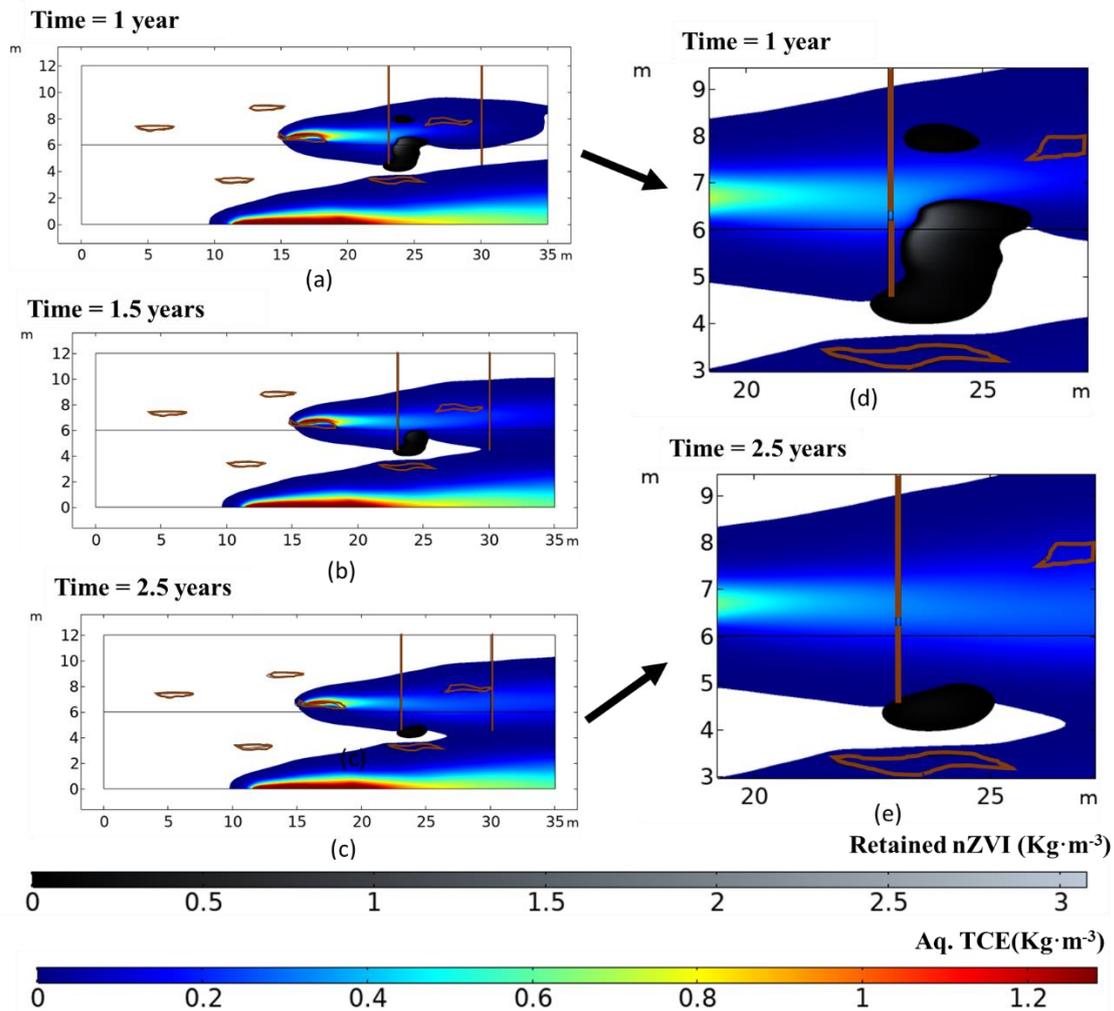

Figure 16: Model prediction for aqueous TCE plume ($kg_{TCE} \cdot m^{-3}$) water and unreacted nZVI plume ($kg_{nZVI} \cdot m^{-3}$) after (a) 5 days (b) 2 months and (c) 6 months since the injection of nZVI. Enlarged version of profile (c) 5 days (d) 6 months since the injection of nZVI, to study the model result at the interface of TCE plume and nZVI mass.

The simulated aqueous TCE concetrnation at the monitoring well is studiesd The monitoring well has the well screen at the depth of 6.6 m, same as the depth of the well screen of the injection well. Figure 17 shows the model predcition for evolution of aqueous TCE at the monitoring well since the nZVI injectiton in the simulated domain. The plot shows that the concetration of simulted aq. TCE dropped to the negegible value after 150 days of injection. The dip in the TCE cocnetration indicates the efficiency of nZVI in TCE degradation resulting in restricted downgradient movement of the TCE plume. The model predicts that the aquesous TCE value in the monitoring well remains neglegible for 1 year since the nZVI injection. Thus the result indicate that the pilot injection in the simulated domain is efficient in the partial removal of TCE for a duration of 1 year. The model predicts the rise of concetration between the first and second year after the injection. The rise in concenreation is 55% less that the aqueous TCE cocnentration at the onset of nZVI injection. This reduction of aqueous TCE cocnetration can be attributed to the depletion of TCE phase acting as a continuous source of aqueous TCE plume. The further decrease in the TCE sourze zone results in the decline of the TCE concentration at the monitoring wel as predicted by the model for the subsequent year.

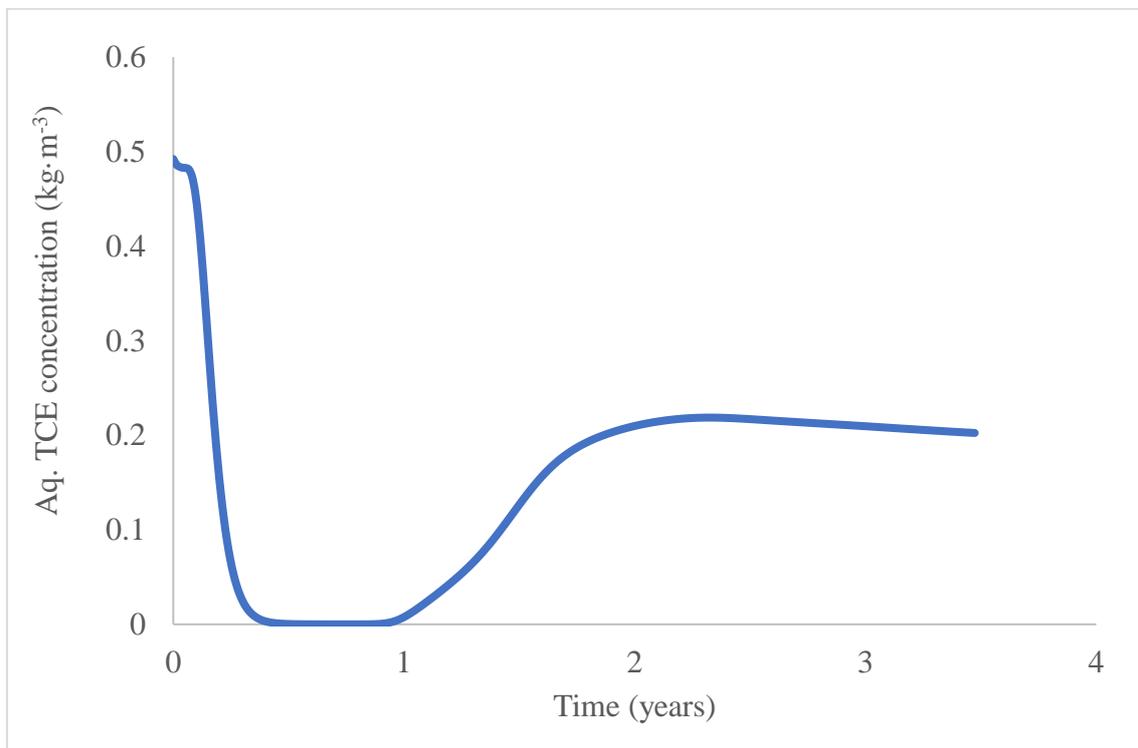

*Figure 17: Model prediction for the aqueous TCE concentration observed at the monitoring well placed at a distance of 7 m from the injection well. The decrease in the TCE concentration followed by the increase indicates the partial degradation of the TCE contamination within the simulated aquifer.*

# 4 Conclusions

In this research work, an integrated framework providing an end-to-end modeling solution has been conceptualized, formulated and deployed for a 2D hypothetical aquifer. The integrated framework can provide accurate and scalable simulation results due to its fully coupled formulations which in addition also enhances its interpretability. The overall objective of the modeling exercise has been to demonstrate the capability of the integrated framework in providing relevant predictive insights at different stages of remediation scheme using nZVI injection. For this purpose, the model uses a 2D hypothetical aquifer, with realistic hydrogeological properties including a natural groundwater flow. The event of contaminant and subsequent planning of remediation have been divided into four stages. The model successfully simulates the primary stage of TCE phase infiltration in the simulated domain. The result predicts TCE source zone architecture which is often challenging to be obtained from field analysis. Using the result of Stage 1 as starting point for Stage 2, the integrated model further simulates the dissolution of TCE and transport of aqueous contaminant aquifer in the model domain (Stage 2). The result demonstrates the simulation's capability in the reconstruction of past groundwater contamination events including its spatio-temporal evolution and depletion of source zone. Furthermore, an essential aspect in a remediation design is the prediction of efficiency of nZVI and timeline of contamination degradation. Therefore, in Stage 3, the transport and retention of nZVI in the model domain has been simulated. The results demonstrate the model prediction for nZVI distribution in the model domain. In addition, the model also characterizes the clogging effect and viscosity changes due to CMC-nZVI injection. Under the integrated framework, the model further estimates the initial impact of CMC-nZVI injection on aqueous TCE plume geometry. In the last step, Stage 4, a reactive transport modeling for aqueous TCE degradation due to nZVI has been developed. The result from this step determines the efficiency of nZVI, and timeline of TCE degradation along with the simulation of reoccurrence of contaminant in the simulated aquifer. The uniqueness of this modeling exercise lies in its integrated approach for the contaminant characterization and remediation design which makes it an end-to-end numerical simulator. Additionally, the model can be upscaled to 3D field scale domain incorporating complex geometry and interaction of groundwater, contaminant, and nZVI. Thus, this research study opens a new numerical approach for the assessment of field sites and development of remediation strategies. Further studies are required to use real field data and calibration for the industrial scale use of the model.